\ifdefined\journalstyle
   \documentclass[review,onefignum,onetabnum]{siamonline250211} 
\else
   \documentclass[11pt, a4paper]{article} 
   
\fi
   
\usepackage{lipsum}
\usepackage{amsfonts}
\usepackage{graphicx}
\usepackage{epstopdf}
\usepackage{algorithmic}
\ifpdf
  \DeclareGraphicsExtensions{.eps,.pdf,.png,.jpg}
\else
  \DeclareGraphicsExtensions{.eps}
\fi

\usepackage{amsmath}
\usepackage{latexsym}
\usepackage{amssymb}
\usepackage{subcaption} 
\usepackage[dvipsnames,svgnames]{xcolor}
\usepackage[acronym]{glossaries}
\newacronym{uq}{UQ}{uncertainty quantification}
\newacronym{mc}{MC}{Monte Carlo}
\newacronym{qmc}{QMC}{quasi-Monte Carlo}
\newacronym{pde}{PDE}{partial differential equation}
\newacronym{mri}{MRI}{magnetic resonance imaging}
\newacronym[firstplural=quantities of interest (QoIs)]{qoi}{QoI}{quantity of interest}
\newacronym{fem}{FEM}{finite element method}
\newacronym{kl}{KL}{Karhunen--Lo{\`e}ve}

\usepackage{enumitem}
\setlist[enumerate]{leftmargin=.5in}
\setlist[itemize]{leftmargin=.5in}

\newtheorem{theorem}{Theorem}[section]

\numberwithin{equation}{section}     
\newenvironment{proof}{\begin{trivlist}
\item[\hskip\labelsep{\it Proof.}]}{$\hfill\Box$\end{trivlist}}
\usepackage{hyperref}
\usepackage{cleveref}

\title{Quasi-Monte Carlo methods for uncertainty quantification of tumor growth modeled by a parametric semi-linear parabolic reaction-diffusion equation\footnote{ADG, FYK, and IHS acknowledge financial support from the Australian Research Council Discovery Project (DP240100769). DN acknowledges funding from Research Foundation Flanders (FWO G091920N). FYK, DN, and IHS acknowledge support from the J. Tinsley Oden Faculty Fellowship Research Program. GP and KW acknowledge support from Department of Energy grant DE-SC002317 and NSF FDT-Biotech award 2436499. The authors acknowledge the Texas Advanced Computing Center (TACC) at The University of Texas at Austin for providing computational resources (Lonestar6) which contributed to the research results reported within this paper.}}

\author{Alexander D.~Gilbert\footnote{School of Mathematics and Statistics, UNSW Sydney, Sydney NSW 2052, Australia 
(\texttt{alexander.gilbert@unsw.edu.au}, \texttt{f.kuo@unsw.edu.au}, \texttt{i.sloan@unsw.edu.au}).}
\and Frances Y.~Kuo\footnotemark[2]
\and Dirk Nuyens\thanks{Department of Computer Science, KU Leuven, Belgium. (\texttt{dirk.nuyens@kuleuven.be}).}
\and Graham Pash\thanks{Oden Institute for Computational Engineering and Sciences, University of Texas at Austin, Austin, TX 78712 USA (\texttt{gtpash@utexas.edu}, \texttt{kwillcox@oden.utexas.edu}).}
\and Ian H.~Sloan\footnotemark[2]
\and Karen E.~Willcox\footnotemark[4]}

\date{July 2026}

\ifpdf
\hypersetup{
pdftitle={Quasi-Monte Carlo methods for uncertainty quantification of tumor growth modeled by a parametric semi-linear parabolic reaction-diffusion equation},
pdfauthor={A.~D.~Gilbert, F.~Y.~Kuo, D.~Nuyens, G.~Pash, I.~H. Sloan, and K.~E.~Willcox}
}
\fi
  
\usepackage{amsopn}

\newcommand{\bsbeta}{{\boldsymbol{\beta}}}
\newcommand{\bsDelta}{{\boldsymbol{\Delta}}}

\newcommand{\bse}{{\boldsymbol{e}}}
\newcommand{\bsrho}{{\boldsymbol{\rho}}}
\newcommand{\bsgamma}{{\boldsymbol{\gamma}}}

\newcommand{\bsnu}{{\boldsymbol{\nu}}}

\newcommand{\bsn}{{\boldsymbol{n}}}

\newcommand{\bsx}{{\boldsymbol{x}}}

\newcommand{\bsy}{{\boldsymbol{y}}}
\newcommand{\bsz}{{\boldsymbol{z}}}

\newcommand{\bsm}{{\boldsymbol{m}}}
\newcommand{\bsk}{{\boldsymbol{k}}}
\newcommand{\bszero}{{\boldsymbol{0}}}

\newcommand{\bshalf}{{\boldsymbol{\tfrac{1}{2}}}}
\newcommand{\rd}{{\mathrm{d}}}

\newcommand{\bbN}{{\mathbb{N}}}

\newcommand{\bbR}{{\mathbb{R}}}
\newcommand{\bbZ}{{\mathbb{Z}}}
\newcommand{\bbE}{{\mathbb{E}}}

\newcommand{\calO}{{\mathcal{O}}}

\newcommand{\calT}{{\mathcal{T}}}

\newcommand{\calX}{{\mathcal{X}}}
\newcommand{\calY}{{\mathcal{Y}}}

\newcommand{\mask}[1]{}
\newcommand{\bhalf}{{[\tfrac{1}{2}]}}

\newcommand{\setu}{{\mathfrak{u}}}

\newcommand{\amax}{a_{\max}}
\newcommand{\amin}{a_{\min}}
\newcommand{\bmax}{b_{\max}}
\newcommand{\bmin}{b_{\min}}
\newcommand{\cmax}{c_{\max}}
\newcommand{\cmin}{c_{\min}}
\newcommand{\kmax}{\kappa_{\max}}
\newcommand{\kmin}{\kappa_{\min}}
\newcommand{\fmax}{f_{\max}}

\definecolor{darkred}{RGB}{139,0,0}
\definecolor{darkgreen}{RGB}{0,100,0}
\definecolor{darkmagenta}{RGB}{180,0,180}
\definecolor{darkblue}{RGB}{0,0,190}
\definecolor{darkorange}{RGB}{180,60,0}
\definecolor{utorange}{RGB}{191,87,0}

\DeclareSymbolFont{bbold}{U}{bbold}{m}{n}
\DeclareSymbolFontAlphabet{\mathbbold}{bbold}

\begin{document}

\maketitle

\begin{abstract}
We study the application of a quasi-Monte Carlo (QMC) method to a class of semi-linear parabolic reaction-diffusion partial differential equations used to model tumor growth. Mathematical models of tumor growth are largely phenomenological in nature, capturing infiltration of the tumor into surrounding healthy tissue, proliferation of the existing tumor, and patient response to therapies, such as chemotherapy and radiotherapy. Considerable inter-patient variability, inherent heterogeneity of the disease, sparse and noisy data collection, and model inadequacy all contribute to significant uncertainty in the model parameters. It is crucial that these uncertainties can be efficiently propagated through the model to compute quantities of interest (QoIs), which in turn may be used to inform clinical decisions. We show that QMC methods can be successful in computing expectations of meaningful QoIs. Well-posedness results are developed for the model and used to show a theoretical error bound for the case of uniform random fields. The theoretical linear error rate, which is superior to that of standard Monte Carlo, is verified numerically. Encouraging computational results are also provided for lognormal random fields, prompting further theoretical development.

\end{abstract}

\ifdefined\journalstyle
  \begin{keywords}
  uncertainty quantification, quasi-Monte Carlo methods, semi-linear parabolic equations, nonlinear analysis, computational oncology
  \end{keywords}
  \begin{MSCcodes}
  65D30, 65D32, 92B05, 92C50, 35K58
  \end{MSCcodes}
\else
  \noindent\textbf{Key words.} uncertainty quantification, quasi-Monte Carlo methods, semi-linear parabolic equations, nonlinear analysis, computational oncology
  
  \medskip
  \noindent\textbf{MSC codes.} 65D30, 65D32, 92B05, 92C50, 35K58
\fi

\section{Introduction}
\Gls*{uq} is crucial in precision medicine, establishing trust in computational models for high-consequence applications. This is particularly salient in oncology, where biomedical imaging plays a critical role in the diagnosis and management of a variety of tumors. However, these data are noisy, infrequently collected, and only indirectly informative, all contributing to significant uncertainty in calibrated models. Practical decision making requires efficient methods for the forward propagation of this uncertainty through mathematical and computational models as well as the estimation of clinically relevant \glspl*{qoi}. 
In this work, we develop \gls*{qmc} methods for efficiently estimating the expectation of \glspl*{qoi} from non-linear \gls*{pde} models of tumor growth.

We consider a semi-linear parabolic reaction-diffusion \gls*{pde} model of tumor growth that captures both invasion into surrounding tissue \cite{alfonso2017biology} and proliferation of the existing lesion \cite{chaplain1996avascular}. This is an established model in the field of computational oncology \cite{swanson2000quantitative} and has been extended to account for patient response to systemic chemo-radiation therapy \cite{hormuth2021image}. Here, we model the uncertain diffusion and proliferation coefficients as random fields, the uncertainty in which leads to uncertainty in the modeled state and \glspl*{qoi}. We will consider the expectation to be integrated over finitely many random variables by truncating a series expansion representing the formally infinite-dimensional input random fields to a sum over finitely many terms. To efficiently address the high-dimensional integral, we will employ \gls*{qmc} methods with the aim of improving upon the slow convergence of standard \gls*{mc} simulation.

\gls*{qmc} methods are a class of quadrature rules designed to compute high-dimensional integrals (or equivalently expectations). The key benefits of QMC methods are that they can be tailored to achieve dimension-independent error bounds with faster than MC convergence rates for smooth integrands. QMC methods have demonstrated impressive performance in several UQ applications involving PDEs with random coefficients or domains, such as linear elliptic \glspl*{pde} \cite{dick2014higher,gantner2018,graham2015quasi,graham2018circulant,guth2021quasi,hakula2024,kuo2012quasi},
elliptic eigenvalue problems \cite{gilbert2019evp,nguyen2024evp}
and more recently linear parabolic equations with control \cite{GKKSS24}.
However, much of the theory has focused on \emph{linear} PDEs. The present work extends the theory in the case of \emph{non-linear}, specifically \emph{semi-linear}, \emph{parabolic} equations, providing a roadmap for adapting \gls*{qmc} methods to wider classes of applications. To the best of our knowledge, this paper is the first to analyze QMC methods for a non-linear PDE.

Well-posedness results are developed for the parametric PDE model in general. Rigorous error estimates are proven for the special case of uniform random fields and confirmed with numerical experiments. 
Note that we do not analyze the error originating from the spatio-temporal discretization of the governing \gls*{pde}, leaving this for future work. Encouraging numerical results are also presented for a more realistic test case of log-normal random fields on a complex, unstructured domain, used previously in models of gliomas in the brain \cite{liang2023bayesian,pash2025predictive}.

One of our key intermediate results in the QMC error analysis is a bound on the derivatives of the solution with respect to uniform random parameters. 
Similar parametric regularity results have also been studied for related classes of semi-linear elliptic PDEs in \cite{CL24b,harbrecht2024,he2024analytic}; see also e.g., \cite{an2025sparse, chernov2025gevrey} for other nonlinear PDE or SPDE.

However, the nonlinear term considered here, namely, $u(1 - u)$ (see \eqref{eq:pde1}), does not fit the settings in \cite{CL24b,harbrecht2024,he2024analytic}. In particular, it is not monotone and so requires a novel proof, which uses the property that if our initial condition lies in $[0, 1]$ then the solution remains in $[0, 1]$ (proven in Section~\ref{sec:apriori}). We also make use of the falling factorial combinatorial bound from~\cite{CL24a}. 
We stress that obtaining regularity bounds for parabolic (time-dependent) nonlinear problems here is significantly harder than for elliptic problems in \cite{CL24a,CL24b,harbrecht2024,he2024analytic} due to the need to handle both space and time derivatives and integrals in the norms.
The parametric regularity bounds presented here are also of independent interest beyond QMC since they are often required to analyze other parametric approximation methods, e.g., sparse grids or stochastic collocation.

The remainder of this paper is organized as follows. In \Cref{sec:formulation} we present the governing equations and variational framework, and then establish well-posedness and an \textit{a priori} bound with explicit constant depending on the parameters for general random field models. A parametric regularity bound is obtained for the case of uniform random fields in \Cref{sec:regularity}. Our main result addresses the \gls*{qmc} error rate 
in \Cref{sec:error}. The theoretical results are confirmed with supporting numerical experiments in \Cref{sec:numerics}. Conclusions are drawn in \Cref{sec:conclusions}.

\section{Problem formulation}
\label{sec:formulation}

\subsection{Governing equation}
\label{sec:pde1}  
The invasion of a solid tumor into surrounding healthy tissue is typically modeled as a diffusion process, while the proliferation of 
existing tumor is modeled with logistic growth subject to some carrying capacity. These phenomenological assumptions give rise to a semi-linear parabolic reaction-diffusion \gls*{pde} \cite{murray2003mathematical}, which has demonstrated predictive value for a variety of solid tumors \cite{clatz2005realistic,hormuth2015predicting,jarrett2021quantitative,lorenzo2019computer}. The model is given by
\begin{align}
\label{eq:pde1}
\mbox{PDE I:}\quad
    \begin{cases}
        \displaystyle\frac{\partial u}{\partial t} 
        - \nabla\cdot (a\,\nabla u) - \kappa\,u\,(1-u) + f\,u
        = 0, \qquad & \bsx\in\Omega,\quad t\in I, \\
        u(\bsx,0) = u_0(\bsx), & \bsx\in\Omega \\
        \nabla u\cdot\bsn = 0, & \bsx\in\partial\Omega,\quad t\in I,
    \end{cases}
\end{align}
where $u=u^\bsy(\bsx,t)$ models the \emph{tumor volume fraction} varying in time $t$ on the interval $I:=[0, T]$ and over spatial coordinates $\bsx\in\Omega$, where $\Omega\subset\mathbb{R}^d$ is the bounded spatial domain with physical dimension $d\in\{2,3\}$ and Lipschitz boundary $\partial\Omega$ with outward unit normal $\bsn$. The homogeneous Neumann boundary condition ensures that the tumor does not grow beyond the boundary of the domain. The linear reaction term $f\,u$ models the effect of treatment. The initial tumor volume fraction $u_0(\bsx)$ may be estimated from medical imaging data. It is natural to assume
\begin{align} \label{eq:pde1-u0}
  0\le u_0(\bsx)\le 1
  \qquad \mbox{for all } \bsx\in\Omega.
\end{align}

Systemic therapies, such as chemotherapy and radiotherapy form the basis of many treatment protocols, for example, gliomas in the brain \cite{stupp2005radiotherapy}. Combination chemoradiation may be modeled as a linear effect on the state, $f = f(\bsx,t) = f_{\rm rt} (\bsx,t) + f_{\rm ct}(\bsx,t)$, with
\begin{align}
  f_{\rm rt} (\bsx, t) &= \gamma_{\text{rt}, \varepsilon}\sum_k 
  \big[1-\exp(-\alpha_{\text{rt}}\, z_{\text{rt},k}(\bsx) - \beta_{\text{rt}}\, z^2_{\text{rt},k}(\bsx))\big]\,\, \mathbbold{1}(t\in[\tau_{\text{rt},k}, \tau_{\text{rt},k} + \varepsilon]), \label{eq:f-rt}
  \\
  f_{\rm ct} (\bsx,t) &= \alpha_{\rm ct}\,\sum_\ell z_{\text{ct},\ell}(\bsx)\, \exp(-\beta_{\rm ct}(t-\tau_{{\rm ct},\ell})) \,\mathbbold{1}(t\geq\tau_{\text{ct},\ell}),
  \label{eq:f-ct}
\end{align}
where $\calT_{\rm rt} = \{\tau_{{\rm rt},k}\}_k$ and $\calT_{\rm ct} = \{\tau_{{\rm ct},\ell}\}_\ell$ denote the collection of times at which radiotherapy and chemotherapy are administered, respectively, and $\mathbbold{1}$ denotes the indicator function. Radiotherapy is assumed to be instantaneous at the moment of treatment with no lasting or time-delaying effects. The application duration is much smaller than the time-discretization, that is, $\varepsilon \ll \Delta t$, where $\Delta t$ is the step size. The surviving fraction is computed with the linear-quadratic model: $z_{{\rm rt},k}(\bsx)$ is the radiotherapy dosage, $\alpha_{\rm rt}>0$ and $\beta_{\rm rt}>0$ are parameters describing the radiosensitivity of the tissue, and $\gamma_{{\rm rt},\varepsilon}>0$ is a parameter that ensures appropriate physical units \cite{borasi2016modelling,mcmahon2018linear,rockne2009mathematical}. Chemotherapy is assumed to follow a decaying exponential model: $z_{{\rm ct},\ell}(\bsx)$ is the drug concentration, $\alpha_{\rm ct}>0$ is the efficacy, and $\beta_{\rm ct}>0$ is the clearance rate \cite{hormuth2021image,jarrett2021quantitative}.

The superscript $\bsy$ denotes a parametric vector modeling any randomness and we allow time-varying coefficients for the diffusion $a = a^\bsy(\bsx,t)$ and proliferation $\kappa = \kappa^\bsy(\bsx,t)$. The precise details on the modeling of the random fields is deferred to \Cref{sec:regularity}. However, we will restrict the parametric vectors $\bsy$ to some admissible parameter set $U \subset\bbR^\bbN$, such that
\begin{align}
\label{eq:pde1-akf}
    \begin{cases}
          0 < \amin^\bsy \le  a^\bsy(\bsx,t)\le \amax^\bsy < \infty, \\
          0 \le\,\kmin^\bsy \le  
          \kappa^\bsy(\bsx,t)\le \kmax^\bsy < \infty, \; \kmax^\bsy > 0, \\
          0 = f_{\min} \le\,  f(\bsx,t)\le \fmax < \infty,
    \end{cases}
    \qquad\mbox{for all } \bsx\in\Omega, \; t\in I, \; \bsy\in U.
\end{align}

\subsection{Functional framework}
\label{sec:function-spaces}
Let $H^1(\Omega)$ denote the usual first-order Sobolev space on $\Omega$ with (squared) norm
\[
 \|v\|_{H^1(\Omega)}^2 := \|v\|_{L^2(\Omega)}^2 + \|\nabla v\|_{L^2(\Omega)}^2.
\]
Let $H^1(\Omega)^*$ denote the dual space of $H^1(\Omega)$, i.e., the space of all bounded linear functionals on $H^1(\Omega)$. We denote the duality pairing between $H^1(\Omega)^*$ and $H^1(\Omega)$ by $\langle \cdot,\cdot\rangle_{H^1(\Omega)^*, H^1(\Omega)}$, with $L^2(\Omega)=L^2(\Omega)^*$ as the pivot space. The dual norm is defined by
\[
  \|z\|_{H^1(\Omega)^*} 
  := \sup_{0\ne v\in H^1(\Omega)} 
  \frac{|\langle z,v \rangle_{H^1(\Omega)^*,H^1(\Omega)}|}{\|v\|_{H^1(\Omega)}}.
\]

It is advantageous to consider functions of time with values in Banach spaces when dealing with time-dependent problems. For a Banach space $B$, let the space $L^2(I;B)$ consist of all measurable functions $v: I\to B$ with finite norm $\|v\|_{L^2(I;B)}^2 := \int_I \|v(\cdot,t)\|_B^2\,\rd t$, and let $C(I;B)$ consist of all functions $v: I\to B$ that are continuous at every $t\in I$ with norm
$\|v\|_{C(I;B)} := \max_{t\in I} \|v(\cdot,t)\|_B$.
Specifically, in this paper we make use of the following norms:
\begin{align*}
  \|v\|_{L^2(I;L^2(\Omega))}^2 
  &:= \int_I \|v(\cdot,t)\|_{L^2(\Omega)}^2\,\rd t, 
  \\
  \|v\|_{L^2(I;H^1(\Omega))}^2 
  &:= \int_I \|v(\cdot,t)\|_{H^1(\Omega)}^2\,\rd t
  = \int_I \big(\|v(\cdot,t)\|_{L^2(\Omega)}^2 + \|\nabla v(\cdot,t)\|_{L^2(\Omega)}^2\big)\,\rd t, 
  \\
  \|v\|_{L^2(I;H^1(\Omega)^*)}^2 
  &:= \int_I \|v(\cdot,t)\|_{H^1(\Omega)^*}^2\,\rd t,
  \\
  \|v\|_{C(I;L^2(\Omega))} 
  &:= \max_{t\in I} \|v(\cdot,t)\|_{L^2(\Omega)}.
\end{align*}
Note that $L^2(I;H^1(\Omega))^* = L^2(I;H^1(\Omega)^*)$. 

We follow Tr\"oltzsch~\cite[\S3.4]{Tro10} to analyze our \gls*{pde} solution in the space $\calX$ (denoted by $W(0,T)$ in \cite{Tro10}) which consists of all functions $v\in L^2(I;H^1(\Omega))$ with (distributional) derivative $\tfrac{\partial}{\partial t} v\in L^2(I;H^1(\Omega)^*)$, and we equip the space $\calX$ with the (graph) norm
\begin{align*}
  \|v\|_\calX^2 
  &:= \|v\|_{L^2(I;H^1(\Omega))}^2 + \|\tfrac{\partial}{\partial t} v\|_{L^2(I;H^1(\Omega)^*)}^2 \\
  &\,= \int_I \big(\|v(\cdot,t)\|_{H^1(\Omega)}^2 
  + \|\tfrac{\partial}{\partial t}v(\cdot,t)\|_{H^1(\Omega)^*}^2\big)\,\rd t \\
  &\,= \int_I \big(\|v(\cdot,t)\|_{L^2(\Omega)}^2 + \|\nabla v(\cdot,t)\|_{L^2(\Omega)}^2 
  + \|\tfrac{\partial}{\partial t}v(\cdot,t)\|_{H^1(\Omega)^*}^2\big)\,\rd t.
\end{align*}

A fundamental property we will use is that the space $\calX$ is embedded in $C(I;L^2(\Omega))$, see \cite[Theorems~3.10 and~3.11]{Tro10}.

We formally write, for $v\in\calX$ and all $\tau\in I$
\begin{align}\label{eq:integration_by_parts}
 \int_0^\tau \!\big\langle \tfrac{\partial}{\partial t} v, v\big\rangle_{H^1(\Omega)^*,H^1(\Omega)}\, \rd t
 = \int_0^\tau \! \frac{1}{2} \frac{\rd}{\rd t} \|v(\cdot,t)\|_{L^2(\Omega)}^2\,\rd t 
 = \frac{\|v(\cdot,\tau)\|_{L^2(\Omega)}^2 
 - \|v(\cdot,0)\|_{L^2(\Omega)}^2}{2}.
\end{align}

We will make use of a number of embeddings which we summarize together here:
\begin{align} 
  \|v\|_{H^1(\Omega)^*} &\le \|v\|_{L^2(\Omega)} \le \|v\|_{H^1(\Omega)},
  \nonumber \\
  \|v\|_{L^4(\Omega)} &\le \theta_1\,\|v\|_{H^1(\Omega)}, \label{eq:emb2} \\
  \|v\|_{C(I;L^2(\Omega))} &\le \theta_2\,\|v\|_\calX, \label{eq:emb3}
\end{align}
with $\theta_1, \theta_2 >0$. The embedding \eqref{eq:emb2} is due to the Sobolev embedding theorem for $H^1(\Omega)$, see e.g., 
\cite[Theorem~7.1]{Tro10}, which states that $H^1(\Omega)$ is embedded in $L^p(\Omega)$ when $d=2$ and $1\le p < \infty$, and when $d\ge 3$ and $1\le p \le \frac{2d}{d-2}$.

Another property we will use is that the space $H^1(\Omega)$ is closed under addition, subtraction, and the maximum and minimum operations, see e.g., \cite[Lemma~7.6]{GT77}.
 
This allows us to define continuous piecewise-smooth test functions in $H^1(\Omega)$.

\subsection{Reparameterization}
\label{sec:pde2}

In the case where $f=0$, that is, untreated growth, the \gls*{pde} model \eqref{eq:pde1} is the well-studied Fisher-KPP equation \cite{fisher1937wave,kolmogoroff1988study}, and well-posedness results are readily available, see e.g., \cite[\S3.10.3]{Per15}. However, we consider $f\neq 0$, thus standard results are not immediately applicable. PDE I \eqref{eq:pde1} can be rewritten as
\[
  \frac{\partial u}{\partial t} - \nabla\cdot (a\, \nabla u) 
  + \kappa\,u^2 + (f-\kappa)\, u = 0.
\]
It is not known \textit{a priori} whether $f\ge \kappa$ holds. Moreover, we are unable to deduce monotonicity of the function $\kappa\,u^2 + (f-\kappa)\, u$ with respect to $u$. Consequently, we cannot claim the existence and uniqueness of a solution using known results, such as those in \cite{Tro10}.

Following a standard technique for treating parabolic \glspl*{pde} involving a linear term with a negative coefficient (see e.g., \cite[\S7.5 Problem~8]{Eva10}), we use the substitution 
\begin{align*} 
  u^\bsy(\bsx,t) = e^{\lambda^\bsy\,t}\, w^\bsy(\bsx,t)\,,
\end{align*}
for some $\lambda^\bsy \in \bbR$, possibly depending on $\bsy$, to obtain
\[
  u = e^{\lambda\,t}\,w \quad\implies\quad
  \frac{\partial u}{\partial t} 
  = \lambda\,e^{\lambda\,t}\,w + e^{\lambda\,t}\,\frac{\partial w}{\partial t}
  \quad\mbox{and}\quad
  \nabla u = e^{\lambda\,t}\,\nabla w,
\]
where to simplify notation we have omitted the $\bsy$ dependence.
Therefore, after dividing through by $e^{\lambda\,t}$ and regrouping, the \gls*{pde} becomes
\[
  \frac{\partial w}{\partial t} 
  - \nabla\cdot (a\,\nabla w) + \kappa\,e^{\lambda\,t}\,w^2 + (\lambda + f-\kappa)\,w = 0.
\]

Thus we have the reparametrized \gls*{pde} problem,
\begin{align}
\label{eq:pde2}
\mbox{PDE II:}\quad
    \begin{cases}
        \displaystyle\frac{\partial w}{\partial t} 
        - \nabla\cdot (a\,\nabla w) + b\,w^2 + c\,w 
        = 0, \qquad & \bsx\in\Omega,\quad t\in I, \\
        w(\bsx,0) = w_0(\bsx), & \bsx\in\Omega \\
        \nabla w\cdot\bsn = 0, & \bsx\in\partial\Omega,\quad t\in I,
    \end{cases}
\end{align}
with coefficients and initial condition related to \gls*{pde} I \eqref{eq:pde1} via
\begin{align} \label{eq:pde2-sub}
  \begin{cases}
  b^\bsy(\bsx,t) := \kappa^\bsy(\bsx,t)\,\exp( \lambda^\bsy\,t),  \\
  c^\bsy(\bsx,t) := \lambda^\bsy + f(\bsx,t) - \kappa^\bsy(\bsx,t), \\
  w_0(\bsx) := u_0(\bsx).
  \end{cases}
\end{align}
To ensure that $c^\bsy(\bsx,t) > 0$, we may take any $\lambda^\bsy > \kmax^\bsy$. From \eqref{eq:pde1-u0} we have $0\le w_0(\bsx)\le 1$ for all $\bsx\in\Omega$. It follows from \eqref{eq:pde1-akf} that
\begin{align} \label{eq:pde2-abc}
  \begin{cases}
  0 < \amin^\bsy \le a^\bsy(\bsx,t) \le \amax^\bsy < \infty, \\
  0 \le \kmin^\bsy \le \bmin^\bsy \le 
  b^\bsy(\bsx,t) \le \bmax^\bsy \le \kmax^\bsy\,e^{\lambda^\bsy\,T} < \infty,  \\
  0 < \lambda^\bsy - \kmax^\bsy \le \cmin^\bsy \le c^\bsy(\bsx,t) \le \cmax^\bsy \le
  \lambda^\bsy + \fmax
  < \infty,
  \end{cases}
 \mbox{for all } \bsx\in\Omega, \; t\in I, \; \bsy\in U.
\end{align}

\subsection{Well-posedness}
\label{sec:well-posed}

In this section, we provide an existence and uniqueness result for the system. We begin by deriving the parametric variational formulation of the PDE problem.

For fixed $\bsy\in U$, a function $u^\bsy \in \calX$ is said to be a weak solution for \gls*{pde}~I~\eqref{eq:pde1} if 
\begin{align} \label{eq:weak1}
  \begin{cases}
  \displaystyle
  \int_I
  \big\langle \tfrac{\partial}{\partial t} u,v\big\rangle_{H^1(\Omega)^*,H^1(\Omega)}
  \,\rd t
  + \!\int_I\!\int_\Omega \!\big( a\,\nabla u\cdot\nabla v
  + (\kappa\,u^2 + (f-\kappa)\,u)\, v \big) \,\rd\bsx\,\rd t = 0
  &\forall\, v\in\calX, \vspace{0.1cm} \\ 
  \displaystyle
  \int_\Omega u(\cdot,0)\,v\,\rd\bsx = \int_\Omega u_0\,v\,\rd\bsx
  & \forall\, v\in L^2(\Omega).
  \end{cases}
\end{align}  
The first condition was obtained by testing the \gls*{pde} with a function $v\in\calX$, integrating over $\bsx\in\Omega$ and $t\in I$, applying integration by parts on the diffusion term, and noting that $\int_{\partial\Omega} a\,\nabla u\cdot\bsn\,v\,\rd\bsx = 0$ due to the Neumann boundary condition. The integral $\int_I\int_\Omega u^2\,v\,\rd\bsx\,\rd t$ is well defined because if $u\in\calX$ then $u^2\in L^2(I; L^2(\Omega))$ via the embedding~\eqref{eq:emb2}.

Analogously, a function $w^\bsy \in \calX$ is said to be a weak solution for \gls*{pde}~II~\eqref{eq:pde2} if 
\begin{align} \label{eq:weak2}
 \begin{cases}
 \displaystyle
  \int_I
  \big\langle \tfrac{\partial}{\partial t} w,v\big\rangle_{H^1(\Omega)^*,H^1(\Omega)}
  \,\rd t
  + \int_I \int_\Omega \big(a\,\nabla w\cdot\nabla v + (b\,w^2+c\,w)\, v\big)\,\rd\bsx\,\rd t = 0
  &\forall\, v\in \calX, \vspace{0.1cm} \\ 
  \displaystyle
  \int_\Omega w(\cdot,0)\,v\,\rd\bsx = \int_\Omega w_0\,v\,\rd\bsx
  &\forall\, v\in L^2(\Omega).
  \end{cases}
\end{align}

\begin{theorem}[existence, uniqueness, boundedness] \label{thm:exist}
For every $\bsy\in U$ and $\lambda^\bsy > \kmax^\bsy > 0$, consider \gls*{pde} II \eqref{eq:pde2} with bounded initial condition $w_0\in L^2(\Omega)$, $0\le w_0(\bsx)\le 1$ for all $\bsx\in\Omega$, and bounded coefficients \eqref{eq:pde2-sub}--\eqref{eq:pde2-abc}.
Then, there exists a unique weak solution $w^\bsy \in \calX$ to \gls*{pde} II. This unique solution is
bounded by $0\le w^\bsy(\bsx,t)\le e^{-\lambda^\bsy\,t}$ for a.e.~$\bsx\in\Omega$ and all $t\in I$, and it satisfies the \textit{a priori} bound $\|w^\bsy\|_{\calX} \le C^\bsy\,\|w_0\|_{L^2(\Omega)}$
with a constant $C^\bsy>0$.

In turn, $u^\bsy = e^{\lambda^\bsy\,t}\,w^\bsy \in\calX$ is the unique weak solution to \gls*{pde} I \eqref{eq:pde1} with bounded initial condition $u_0 = w_0\in L^2(\Omega)$, $0\le u_0(\bsx)\le 1$ for all $\bsx\in\Omega$, and bounded coefficients \eqref{eq:pde1-akf}, and it is bounded by $0\le u^\bsy(\bsx,t)\le 1$ for a.e.~$\bsx\in\Omega$ and all $t\in I$.
\end{theorem}

Before proving this theorem, we first remark that there is considerable theory developed for the case of monotone semi-linear equations \cite{Tro10}, which we could readily make use of were the term $b\,w^2 + c\,w$ of \gls*{pde} II~\eqref{eq:pde2} monotone in $w$. Unfortunately, $b\,w^2 + c\,w$ is generally not monotone, unless it is known that $w\ge 0$. Moreover, recall that the solution for \gls*{pde} I~\eqref{eq:pde1} represents the tumor volume fraction and should naturally satisfy $0\le u\le 1$, corresponding to $0\le w \le e^{-\lambda\,t}$. But these bounds are not immediately obvious.
We therefore modify our \gls*{pde} II by imposing a restriction on the term $b\,w^2+c\,w$ whenever $w<0$ or $w> e^{-\lambda\,t}$:
\begin{align} \label{eq:pde3}
 \mbox{\gls*{pde} III:}\quad
 \begin{cases}
  \displaystyle\frac{\partial w}{\partial t} 
  - \nabla\cdot (a\,\nabla w) + d(w)
  = 0, \qquad & \bsx\in\Omega,\quad t\in I, \\
  w(\bsx,0) = w_0(\bsx), & \bsx\in\Omega \\
  \nabla w\cdot\bsn = 0, & \bsx\in\partial\Omega,\quad t\in I,
  \end{cases}
\end{align}
where we define a continuous piecewise cutoff function
\begin{align} \label{eq:d-def}
  d(w) :=
  \begin{cases}
  0 & \mbox{if } w<0, \\
  b\,w^2 + c\,w & \mbox{if } 0\le w\le e^{-\lambda\,t}, \\
  b\,e^{-2\lambda\,t} + c\,e^{-\lambda\,t}
  = (\lambda + f)\,e^{-\lambda\,t} 
  & \mbox{if } w>e^{-\lambda\,t}.
  \end{cases}
\end{align}
A weak solution $w^\bsy\in\calX$ to \gls*{pde} III \eqref{eq:pde3} satisfies \eqref{eq:weak2} with $b\,w^2 + c\,w$ replaced by $d(w)$.

Our argument for proving Theorem~\ref{thm:exist} goes as follows: 
\begin{enumerate}
\item For fixed $\bsy\in U$ and $\lambda^\bsy > \kmax^\bsy >0$, we infer from 
\cite[Assumptions~5.1 and~5.2, Lemmas~5.3 and~7.10]{Tro10} that there exists a unique weak solution $w_*^\bsy\in\calX$ of \gls*{pde} III with bounded coefficients \eqref{eq:pde2-abc}, with \textit{a priori} bound $\|w_*^\bsy\|_{\calX} \le C^\bsy\,\|w_0\|_{L^2(\Omega)}$ for $C^\bsy>0$. 
\item We prove that if the initial condition satisfies $0\le w_0(\bsx)\le 1$ for all $\bsx\in\Omega$ then this unique solution $w_*^\bsy\in\calX$
of \gls*{pde}~III satisfies $0\le w_*^\bsy(\bsx,t)\le e^{-\lambda^\bsy\,t}$ for a.e.~$\bsx\in\Omega$ and all $t\in I$.

\item Then, for this unique solution $w_*^\bsy\in\calX$, \gls*{pde}~III holds with $d(w)$ defined only by the middle case. This shows that $w_*^\bsy\in\calX$ is also a solution for \gls*{pde}~II. This confirms the existence of a solution for \gls*{pde}~II satisfying $0\le w_*^\bsy(\bsx,t)\le e^{-\lambda^\bsy\,t}$ for a.e.~$\bsx\in\Omega$ and all $t\in I$, given the bounded initial condition $0\le w_0(\bsx)\le 1$ for all $\bsx\in\Omega$.

\item Finally we prove that this solution $w_*^\bsy\in\calX$ for \gls*{pde}~II is unique. In particular, there is no other solution taking values below $0$ or above $e^{-\lambda^\bsy\,t}$. Hence we conclude that $u^\bsy(\bsx,t) = e^{\lambda^\bsy\,t}\,w_*^\bsy(\bsx,t)$ is the unique solution of \gls*{pde} I, and it satisfies $0\le u(\bsx,t)\le 1$ for a.e.~$\bsx\in\Omega$ and all $t\in I$.
\end{enumerate}

The strategy of working with a cutoff function in \gls*{pde} III and using continuous piecewise-smooth test functions (see our proof below) appears already in the literature, sometimes referred to as Stampacchia truncation,  see e.g., \cite{CCMR25,CGLMRR20,Per15,Tro10}. However, the problem settings in those references differ from ours, so we include a proof for completeness. 

\begin{proof}[Proof of Theorem~\ref{thm:exist}] We begin by proving the existence of a solution.

\textsc{Existence.}
We fix $\bsy\in U$ and assume bounded coefficients \eqref{eq:pde2-abc}.
We write the cutoff function \eqref{eq:d-def} showing explicit dependence on all arguments,
\begin{align*} 
  d(w) = d(\bsx,t,w) :=
  \begin{cases}
  0 & \mbox{if } w^\bsy(x,t)<0, \\
  b^\bsy(\bsx,t)\,\big(w^\bsy(\bsx,t)\big)^2 + c^\bsy(\bsx,t)\,w^\bsy(x,t) 
  & \mbox{if } 0\le w^\bsy(x,t)\le  e^{-\lambda^\bsy\,t}, \\
  b^\bsy(\bsx,t)\, e^{-2\lambda^\bsy\,t}
  + c^\bsy(\bsx,t)\, e^{-\lambda^\bsy\,t} 
  & \mbox{if } w^\bsy(x,t)> e^{-\lambda^\bsy\,t}.
  \end{cases}
\end{align*}

We see that \gls*{pde} III satisfies \cite[Assumption~5.1]{Tro10}:
$\Omega$ is a bounded Lipschitz domain;
for fixed $w\in\bbR$, the function $d(\bsx,t,w)$ is measurable 
with respect to $\bsx\in\Omega$ and $t\in I$;
for all $\bsx\in\Omega$ and $t\in I$, the function $d(\bsx,t,w)$ is monotonically increasing (nondecreasing) with respect to $w\in\bbR$.
We see that \gls*{pde} III also satisfies \cite[Assumption~5.2]{Tro10}:
for all $\bsx\in\Omega$ and $t\in I$, the function $d(\bsx,t,w)$ is uniformly bounded and Lipschitz in $w\in\bbR$.

Hence, by \cite[Lemma~5.3]{Tro10} and subsequent discussions, there exists a unique weak solution of \gls*{pde}~III in $\calX$. Consequently, by \cite[Lemma~7.10]{Tro10}, this solution satisfies 
$\|w^\bsy\|_{\calX} \le C^\bsy\,\|w_0\|_{L^2(\Omega)}$
with a constant $C^\bsy>0$.

We proceed to prove that this unique weak solution for \gls*{pde}~III is bounded between $0$ and~$e^{-\lambda\,t}$ if the initial condition is bounded by $0\le w_0(\bsx)\le 1$ for all $\bsx\in\Omega$. Integrating \eqref{eq:pde3} with respect to $\bsx\in\Omega$, using integration by parts on the diffusion term and using the zero Neumann boundary condition, and then integrating over $t\in [0,\tau]$ for some $\tau\le T$, we obtain
\begin{align} \label{eq:start}
  \int_0^\tau\! \big\langle \tfrac{\partial}{\partial t} w,v\big\rangle_{H^1(\Omega)^*,H^1(\Omega)}\,\rd t
  + \int_0^\tau\! \int_\Omega \big( a\,\nabla w\cdot\nabla v + d(w)\, v \big) \,\rd\bsx\,\rd t = 0
  \quad\forall\, v\in \calX.
\end{align}%

\textsc{Lower bound.}
We define $v = w_\# \in \calX$ to be the ``unwanted lower part'' of $w$,
\begin{align*}
  w_\#(\bsx,t) := \min(w(\bsx,t),0) = 
  \begin{cases}
  w(\bsx,t) & \mbox{if } w(\bsx,t) < 0, \\
  0 & \mbox{otherwise}.
  \end{cases}
\end{align*}
For a.e.~$\bsx\in\Omega$, since $w(\bsx,0) = w_0(\bsx)\ge 0$, we have $w_\#(\bsx,0) = 0$ and $\|w_\#(\cdot,0)\|_{L^2(\Omega)} = 0$. 
If $w(\bsx,t)< 0$ then
\begin{align} \label{eq:simp0}
  \tfrac{\partial}{\partial t} w = \tfrac{\partial}{\partial t} w_\#, \quad
  \nabla w\cdot \nabla w_\# = |\nabla w_\#|^2, \quad
  d(w)\, w_\# = 0,
\end{align}
where we used Case~1 of \eqref{eq:d-def}. On the other hand, if $w(\bsx,t)\ge 0$ then $w_\#(\bsx,t)=0$ and the equalities in \eqref{eq:simp0} hold trivially.
Putting \eqref{eq:simp0} into \eqref{eq:start} and using \eqref{eq:integration_by_parts} gives for all $\tau\in I$,
\begin{align*} 
  &\frac{\|w_\#(\cdot,\tau)\|_{L^2(\Omega)}^2 - \|w_\#(\cdot,0)\|_{L^2(\Omega)}^2 }{2} \\
  &= \int_0^\tau\! \frac{1}{2} \frac{\rd}{\rd t} \|w_\#(\cdot,t)\|_{L^2(\Omega)}^2\,\rd t 
  = \int_0^\tau\! 
  \big\langle \tfrac{\partial}{\partial t} w_\#,w_\# \big\rangle_{H^1(\Omega)^*,H^1(\Omega)} \,\rd t
  = - \int_0^\tau\!\int_\Omega a\,|\nabla w_\#|^2\,\rd\bsx \,\rd t
  \le 0.
\end{align*}%
Thus $\|w_\#(\cdot,\tau)\|_{L^2(\Omega)} \le \|w_\#(\cdot,0)\|_{L^2(\Omega)} = 0$, and hence
$w_\#(\bsx,\tau) = 0$ for a.e.~$\bsx\in\Omega$, and in turn this yields
\begin{align*} 
  w(\bsx,\tau) \ge 0 \qquad\mbox{for a.e.\ } \bsx\in\Omega
  \mbox{ and all } \tau\in I.
\end{align*}  

\textsc{Upper bound.}
Now we define $v = w_\# \in \calX$ on the ``unwanted upper part'' of~$w$,
\begin{align*}
  w_\#(\bsx,t) := \max(w(\bsx,t)-e^{-\lambda\,t},0) =
  \begin{cases}
  w(\bsx,t)- e^{-\lambda\,t} 
  & \mbox{if } w(\bsx,t) >  e^{-\lambda\,t}, \\
  0 & \mbox{otherwise}.
  \end{cases}
\end{align*}
For a.e.~$\bsx\in\Omega$, since  $w(\bsx,0) = w_0(\bsx)\le 1$, we have $w_\#(\bsx,0) = 0$ and $\|w_\#(\cdot,0)\|_{L^2(\Omega)} = 0$.
If $w(\bsx,t)> e^{-\lambda\,t}$ then
$\frac{\partial}{\partial t}  w_\#= \frac{\partial}{\partial t} w + \lambda\,e^{-\lambda\,t}$, and so
\begin{align} \label{eq:simp1}
  \tfrac{\partial}{\partial t}w = \tfrac{\partial}{\partial t} w_\# - \lambda\,e^{-\lambda\,t}, \quad
  \nabla w\cdot \nabla w_\# = |\nabla w_\#|^2, \quad
  d(w)\, w_\# = (\lambda+f)\,e^{-\lambda\,t}\,w_\#,
\end{align}
where we used Case~3 of \eqref{eq:d-def}. On the other hand, if $w(\bsx,t)\le e^{-\lambda\,t}$ then $w_\#(\bsx,t)=0$ and the equalities in \eqref{eq:simp1} hold trivially.
Putting \eqref{eq:simp1} into \eqref{eq:start} gives for all $\tau\in I$,
\begin{align*} 
  \int_0^\tau\!\! \big\langle \tfrac{\partial}{\partial t} w_\# - \lambda\,e^{-\lambda\,t},  w_\#\big\rangle_{H^1(\Omega)^*,H^1(\Omega)}\,\rd t
  + \int_0^\tau\!\! \int_\Omega \big( a\,|\nabla w_\#|^2 + (\lambda+f)\,e^{-\lambda\,t}\, w_\#\big)
  \,\rd\bsx\,\rd t = 0.
\end{align*}%
Using \eqref{eq:integration_by_parts} then yields
\begin{align*} 
  \frac{\|w_\#(\cdot,\tau)\|_{L^2(\Omega)}^2 - \|w_\#(\cdot,0)\|_{L^2(\Omega)}^2 }{2}
  = - \int_0^\tau\!\! \int_\Omega \big(a\,|\nabla w_\#|^2 + f\,e^{-\lambda\,t} \,w_\#\big)
  \,\rd\bsx\,\rd t
  \le 0.
\end{align*}%
Thus $\|w_\#(\cdot,\tau)\|_{L^2(\Omega)}  \le \|w_\#(\cdot,0)\|_{L^2(\Omega)} = 0$, and hence
$w_\#(\bsx,\tau) = 0$ for a.e.~$\bsx\in\Omega$, and in turn this yields
\begin{align*} 
  w(\bsx,\tau) \le e^{-\lambda\,t} \qquad\mbox{for a.e.\ } \bsx\in\Omega
  \mbox{ and all } \tau\in I.
\end{align*}  

\textsc{Uniqueness.}
We have shown that the unique solution of \gls*{pde} III satisfies $0\le w(\bsx,t) \le e^{-\lambda\,t}$ when the initial condition satisfies $0\le w_0(\bsx) \le 1$. For this solution, \gls*{pde} III holds with $d(w)$ defined only by Case~2 in \eqref{eq:d-def}. Thus, this solution for \gls*{pde} III is also a solution for \gls*{pde} II, demonstrating the existence of a solution for \gls*{pde} II satisfying $0\le w(\bsx,t) \le e^{-\lambda\,t}$. 

Suppose that $w_2$ is another solution for \gls*{pde} II, which is not necessarily bounded between $0$ and $e^{-\lambda\,t}$. Subtracting their respective equations analogous to \eqref{eq:start}, we obtain for all $\tau\in I$,
\begin{align*} 
  &\int_0^\tau\! \big\langle \tfrac{\partial}{\partial t} (w-w_2),v\big\rangle_{H^1(\Omega)^*,H^1(\Omega)}\,\rd t \\
  &+ \int_0^\tau\! \int_\Omega \big( a\,\nabla (w-w_2)\cdot\nabla v + b\,(w^2-w_2^2)\,v + c\,(w-w_2)\, v \big) \,\rd\bsx\,\rd t = 0
  \quad\forall\, v\in \calX.
\end{align*}
Taking $v = w - w_2$ and writing $w^2 - w_2^2 = 2\,w\,(w-w_2) - (w-w_2)^2 = 2\,w\,v - v^2$, we obtain
\begin{align*} 
  \int_0^\tau\! \big\langle \tfrac{\partial}{\partial t} v,v\big\rangle_{H^1(\Omega)^*,H^1(\Omega)}\,\rd t + \int_0^\tau\!\! \int_\Omega \big( a\,|\nabla v|^2 + (2\,b\,w + c)\,v^2 - b\,v^3 \big) \,\rd\bsx\,\rd t = 0.
\end{align*}
Since $\|v(\cdot,0)\|_{L^2(\Omega)}^2 = 0$ and $w \ge 0$, and using \eqref{eq:integration_by_parts}, we conclude that
\begin{align} \label{eq:mango} 
  \frac{\|v(\cdot,\tau)\|_{L^2(\Omega)}^2}{2}
 + \amin\, \|\nabla v\|_{L^2([0,\tau];L^2(\Omega))}^2 + \cmin\, \|v\|_{L^2([0,\tau];L^2(\Omega))}^2 
 \le b_{\max} \int_0^\tau\!\! \int_\Omega \vert v\vert^3 \,\rd\bsx\,\rd t.
\end{align}

We now prove that $v\equiv 0$ using a proof by contradiction. To facilitate this proof, we introduce the space $\calY_\tau := C([0,\tau];L^2(\Omega)) \cap L^2([0,\tau];H^1(\Omega))$ for $\tau\in [0,T]$ with the norm
\begin{align*}
  \|v\|_{\calY_\tau}^2 :=
  \max_{t\in [0,\tau]} \|v(\cdot,t)\|_{L^2(\Omega)}^2 
  + \|\nabla v\|_{L^2([0,\tau];L^2(\Omega))}^2 + \|v\|_{L^2([0,\tau];L^2(\Omega))}^2.
\end{align*}
Let $t_0 := \max\{\tau\in [0,T] : \|v\|_{\calY_\tau} = 0\}$. Then $v(\cdot,t) = w(\cdot,t) - w_2(\cdot,t) = 0$ for all $t\le t_0$. We assume by contradiction that $t_0<T$. For any $\varepsilon>0$, we take $t_1\in (t_0,T]$ such that $0 < \|v\|_{\calY_{t_1}} < \varepsilon$.

For any $\tau \le t_1$, the right-hand side of \eqref{eq:mango} can be bounded by $\bmax\,\theta_1^2\,\theta_2\,\|v\|_{\calY_{t_1}}^3$ using \eqref{eq:fantastic}. Using this to bound each of the three terms on the left-hand side of \eqref{eq:mango} separately and then combining, we deduce that
\begin{align*}
  \|v\|_{\calY_{t_1}}^2 \le \big(2 + \tfrac{1}{\amin} + \tfrac{1}{\cmin}\big)\, 
  \bmax\,\theta_1^2\,\theta_2\,\|v\|_{\calY_{t_1}}^3
  < \big(2 + \tfrac{1}{\amin} + \tfrac{1}{\cmin}\big)\, 
  \bmax\,\theta_1^2\,\theta_2\,\,\varepsilon\,\|v\|_{\calY_{t_1}}^2,
\end{align*}
which yields a contradiction if we take $\varepsilon = [(2 + \tfrac{1}{\amin} + \tfrac{1}{\cmin})\, 
  \bmax\,\theta_1^2\,\theta_2]^{-1}$. Hence $t_0 = T$, that is, $w \equiv w_2$ and therefore the solution is unique.
\end{proof}

\subsection{A priori bound}
\label{sec:apriori}

Here we derive \textit{a priori} bounds with explicit constants. 

\begin{theorem}[\textit{a priori} bound] \label{thm:priori}
For every $\bsy\in U$, let $u^\bsy\in\calX$ denote the unique weak solution to \gls*{pde} I \eqref{eq:pde1} with bounded initial condition $u_0\in L^2(\Omega)$, $0\le u_0(\bsx)\le 1$ for all $\bsx\in\Omega$, and bounded coefficients \eqref{eq:pde1-akf}. Then we have the \textit{a priori} bound
\begin{align*}
\|u^\bsy\|_{\calX} &\le C^\bsy\,\|u_0\|_{L^2(\Omega)},
 \quad\mbox{with}\quad
 C^\bsy
 := \frac{(1+ \amax^\bsy + \kmax^\bsy + \fmax)\,e^{\kmax^\bsy T}}{\sqrt{\min(\amin^\bsy,\kmax^\bsy)}}.
\end{align*}
\end{theorem}

\begin{proof}
From \gls*{pde}~I \eqref{eq:pde1}, we have for a.e.~$t\in I$ and all $v\in H^1(\Omega)$,
\begin{align} \label{eq:again}
  \big\langle \tfrac{\partial}{\partial t} u,v\big\rangle_{H^1(\Omega)^*,H^1(\Omega)}
  + \int_\Omega \big(a\,\nabla u\cdot\nabla v + (\kappa\,u^2+ (f-\kappa)\,u)\, v\big)\,\rd\bsx = 0.
\end{align} 

\textsc{Time derivative.}
Since $0\le u\le 1$, we have $|(\kappa\,u^2+ (f-\kappa)\,u)\, v| 
= |\kappa\,u\,(u-1)+ f\,u|\, |v|
\le (\kappa+f)\,u\,|v|$. Thus we obtain from \eqref{eq:again}
\begin{align*} 
  \big|\big\langle \tfrac{\partial}{\partial t} u,v\big\rangle_{H^1(\Omega)^*,H^1(\Omega)}\big|
  &\le \amax\, \|\nabla u(\cdot,t)\|_{L^2(\Omega)}\, \|\nabla v\|_{L^2(\Omega)} 
   + (\kmax + \fmax)\,\|u(\cdot,t)\|_{L^2(\Omega)}\,\|v\|_{L^2(\Omega)} \\
  &\le (\amax + \kmax + \fmax )\,\|u(\cdot,t)\|_{H^1(\Omega)}\,\|v\|_{H^1(\Omega)},
\end{align*}  
and therefore
\begin{align*}
  \|\tfrac{\partial u}{\partial t}(\cdot,t)\|_{H^1(\Omega)^*} 
  = \sup_{0\ne v\in H^1(\Omega)} 
  \frac{|\langle \tfrac{\partial}{\partial t} u,v\rangle_{H^1(\Omega)^*,H^1(\Omega)}|}{\|v\|_{H^1(\Omega)}} 
  \le (\amax + \kmax + \fmax )\,\|u(\cdot,t)\|_{H^1(\Omega)},
\end{align*}
which yields
\begin{align} \label{eq:magic1}
  \|\tfrac{\partial}{\partial t}u\|_{L^2(I;H^1(\Omega)^*)} 
  \le (\amax + \kmax + \fmax  )\,\|u\|_{L^2(I;H^1(\Omega))}.
\end{align}

\textsc{Gradient.}
Taking $v = u(\cdot,t)$ in \eqref{eq:again} and integrating with respect to $t\in [0,\tau]$ for some $\tau\le T$, and using \eqref{eq:integration_by_parts}, we obtain
\begin{align*} 
  &\frac{\|u(\cdot,\tau)\|_{L^2(\Omega)}^2 - \|u(\cdot,0)\|_{L^2(\Omega)}^2}{2}
  + \int_0^\tau\!\! \int_\Omega \big(
  a\,|\nabla u|^2 + \kappa\,u^2+ (f-\kappa)\,u\big)\,u\,\rd\bsx\,\rd t = 0, 
\end{align*} 
which rearranges into
\begin{align} \label{eq:banana}
  &\frac{\|u(\cdot,\tau)\|_{L^2(\Omega)}^2}{2} 
  + \!\int_0^\tau\!\! \int_\Omega \big (a\,|\nabla u|^2 + f\,u^2 \big) \,\rd\bsx\,\rd t 
  = \!\int_0^\tau\!\! \int_\Omega \kappa\,u^2\,(1-u)\,\rd\bsx\,\rd t
  + \frac{\|u_0\|_{L^2(\Omega)}^2}{2}. 
\end{align}
Dropping the nonnegative second term in \eqref{eq:banana} and using $0\le u\le 1$, we obtain the inequality
\begin{align*} 
  \|u(\cdot,\tau)\|_{L^2(\Omega)}^2
  \le 2\,\kmax \int_0^\tau \|u(\cdot,t)\|_{L^2(\Omega)}^2\,\rd t
  + \|u_0\|_{L^2(\Omega)}^2.
\end{align*} 
The integral form of Gr\"onwall's inequality (see e.g., \cite[Lemma 1.4.1]{quarteroni1994numerical}) then yields
\begin{align*} 
  \|u(\cdot,\tau)\|_{L^2(\Omega)}^2
  \le e^{2\,\kmax\,\tau} \|u_0\|_{L^2(\Omega)}^2,
  \quad\mbox{and so}\quad
  \|u\|_{L^2(I;L^2(\Omega))}^2
  \le \frac{e^{2\,\kmax\,T}-1}{2\,\kmax} \|u_0\|_{L^2(\Omega)}^2.
\end{align*} 
Keeping only the $a\,|\nabla u|^2$ term on the left-hand side of \eqref{eq:banana} and taking $\tau=T$ lead to
\begin{align*} 
  \amin\, \|\nabla u\|_{L^2(I;L^2(\Omega))}^2
  &\le \kmax\, \|u\|_{L^2(I;L^2(\Omega))}^2
  + \frac{\|u_0\|_{L^2(\Omega)}^2}{2}
  \le \frac{e^{2\,\kmax\,T}}{2}\, \|u_0\|_{L^2(\Omega)}^2.
\end{align*} 
In turn we obtain
\begin{align} \label{eq:magic2}
  &\|u\|_{L^2(I;H^1(\Omega))}^2 
  = \|u\|_{L^2(I;L^2(\Omega))}^2 +  \|\nabla u\|_{L^2(I;L^2(\Omega))}^2 \nonumber\\
  &\le \frac{e^{2\,\kmax\,T}-1}{2\,\kmax} \|u_0\|_{L^2(\Omega)}^2
  + \frac{e^{2\,\kmax\,T}}{2\,\amin}\, \|u_0\|_{L^2(\Omega)}^2 
  \le \frac{e^{2\,\kmax\,T}}{\min(\amin,\kmax)}
  \|u_0\|_{L^2(\Omega)}^2.
\end{align}

\textsc{a priori bound.}
Combining \eqref{eq:magic1} and \eqref{eq:magic2}, we obtain
\begin{align*}
  &\|u\|_{\calX}
  = \sqrt{\|u\|_{L^2(I;H^1(\Omega))}^2 + \|\tfrac{\partial}{\partial t}u\|_{L^2(I;H^1(\Omega)^*)}^2} 
  \le \|u\|_{L^2(I;H^1(\Omega))} + \|\tfrac{\partial}{\partial t}u\|_{L^2(I;H^1(\Omega)^*)}\\
  &\le (1 + \amax + \kmax + \fmax)\,\|u\|_{L^2(I;H^1(\Omega))} 
  \le \frac{(1+ \amax + \kmax + \fmax)\, 
  e^{\kmax\,T}}{\sqrt{\min(\amin,\kmax)}} \,\|u_0\|_{L^2(\Omega)},
\end{align*}
as required.
\end{proof}

\section{Parametric regularity}
\label{sec:regularity}

\subsection{Uniform random field model}
We model the coefficients as affine series expansions in terms of uniform parameters $\bsy\in [-\tfrac{1}{2},\tfrac{1}{2}]^\bbN$, see e.g., \cite{GKKSS24,kuo2012quasi},
\begin{align} 
 a^\bsy(\bsx,t) 
 &:= \psi_0(\bsx,t) + \sum_{j\ge 1} y_j\, \psi_j(\bsx,t), \label{eq:a-def}
 \\
 \kappa^\bsy(\bsx,t) 
 &:= \xi_0(\bsx,t) + \sum_{j\ge 1} y_j\, \xi_j(\bsx,t), \label{eq:k-def}
\end{align}
and both of them are bounded independently of $\bsy$,
\begin{align*}
  &0 < \amin := \min_{t\in I} \inf_{\bsx\in\Omega}|\psi_0(\bsx,t)|
  - \frac{1}{2} \sum_{j\ge 1} \|\psi_j\|_\infty
  < \amax := \|\psi_0\|_\infty + \frac{1}{2} \sum_{j\ge 1} \|\psi_j\|_\infty < \infty, \\
  &0 \le \kmin := \min_{t\in I} \inf_{\bsx\in\Omega}|\xi_0(\bsx,t)|
  - \frac{1}{2} \sum_{j\ge 1} \|\xi_j\|_\infty
  < \kmax := \|\xi_0\|_\infty + \frac{1}{2} \sum_{j\ge 1} \|\xi_j\|_\infty < \infty, 
\end{align*}
where we define
\[
  \|v\|_\infty := \max_{t\in I} \|v(\cdot,t)\|_{L^\infty(\Omega)}.
\]
To make the two random fields completely independent, we {will later define $\psi_j \equiv 0$ when $j$ is even and $\xi_j \equiv 0$ when $j$ is odd. This can be interpreted as interleaving two series expansions.

\subsection{Parametric regularity for uniform random fields}
We will write derivatives with respect to the parameters
$\bsy$ using multiindex notation, where for a multiindex $\bsnu \in \bbN^\bbN_0$ with $|\bsnu| := \sum_{j \geq 1} \nu_j < \infty$, we define
the derivative of order $\bsnu$ by
\[
\partial^\bsnu_\bsy = \prod_{j \geq 1} \Big(\frac{\partial}{\partial y_j}\Big)^{\nu_j}.
\]
We will use the \emph{absolute value of the falling factorial} as introduced in Chernov and L\^e~\cite{CL24a,CL24b} 
\[
  \bhalf_0 := 1, \qquad
  \bhalf_n := \big| \tfrac{1}{2}\, (\tfrac{1}{2}-1)\, (\tfrac{1}{2}-2)\cdots (\tfrac{1}{2}-n+1) \big|
  \quad\mbox{for } n\ge 1,
\] 
with $\bhalf_n \le n! \le 2\cdot 2^n\,\bhalf_n$, and we will use a key property from \cite[Lemma~2.3]{CL24b},
\begin{align} \label{eq:falling}
 \sum_{\bszero < \bsm <\bsnu} \binom{\bsnu}{\bsm}\, \bhalf_{|\bsm|}\, \bhalf_{|\bsnu-\bsm|}
 \le 2\,\bhalf_{|\bsnu|},
\end{align}
where $\bsm <\bsnu$ means that $\bsm\ne\bsnu$ and $m_j\le \nu_j$ for all $j\ge 1$, and $\binom{\bsnu}{\bsm} = \prod_{j\ge 1} \binom{\nu_j}{m_j}$ is a product of binomial coefficients.

\begin{theorem}[parametric regularity for uniform model] \label{thm:reg}
For every $\bsy\in U$, let $u^\bsy\in\calX$ denote the unique weak solution to \gls*{pde} I \eqref{eq:pde1} with bounded initial condition $u_0\in L^2(\Omega)$, $0\le u_0(\bsx)\le 1$ for all $\bsx\in\Omega$, and bounded coefficients \eqref{eq:pde1-akf} defined by affine expansions \eqref{eq:a-def} and \eqref{eq:k-def}.
Then for every multiindex $\bsnu$, the $\bsnu$-th derivative with respect to $\bsy$ of the solution satisfies $\partial_\bsy^\bsnu u^\bsy \in \calX$ with the regularity bound
\begin{align*}
 \|\partial^\bsnu_\bsy u^\bsy\|_{\calX} 
 \le \widetilde{C}\,(\rho\,\bsbeta)^\bsnu\,\bhalf_{|\bsnu|}\,\|u_0\|_{L^2(\Omega)}, \quad
  \widetilde{C} := \max(C,\tfrac{E}{\rho}), 
\end{align*}
with the sequence
\begin{align*} 
\beta_j := \max(\|\psi_j\|_\infty,\|\xi_j\|_\infty),
\end{align*}
where $C$ is defined in Theorem~\ref{thm:priori} (in the uniform model, $\amax$, $\amin$, $\kmax$, $\fmax$ are independent of $\bsy$), while
$E$ and $\rho$ are defined in \eqref{eq:DE} and \eqref{eq:rho} below, respectively, both depending on $\Phi$ defined in \eqref{eq:Phi} below.
\end{theorem}

\begin{proof}
For any multiindex~$\bsnu\ne\bszero$, we have the mixed derivative $\partial^\bsnu_\bsy$ of \eqref{eq:a-def} and \eqref{eq:k-def}
\begin{align} \label{eq:diff-ak}
 \partial^\bsnu_\bsy a^\bsy(\bsx,t) 
 = \begin{cases}
 \psi_j(\bsx,t) & \mbox{if $\bsnu = \bse_j$}, \\
 0 & \mbox{otherwise},
 \end{cases} 
 \quad\mbox{and}\quad
 \partial^\bsnu_\bsy \kappa^\bsy(\bsx,t) 
 = \begin{cases}
 \xi_j(\bsx,t) & \mbox{if $\bsnu = \bse_j$}, \\
 0 & \mbox{otherwise}.
 \end{cases}
\end{align}

Recall that the weak solution $u\in\calX$ satisfies \eqref{eq:weak1}. For any multiindex $\bsnu\ne\bszero$ we take the mixed derivative $\partial_\bsy^\bsnu$ on both sides of the first condition in \eqref{eq:weak1}, but we integrate over $t\in [0,\tau]$ for some $\tau\le T$ instead of the full interval $I$. Using the Leibniz product rule, the mixed derivative $\partial_\bsy^\bsnu$ for the space-time integrand in \eqref{eq:weak1} can be written as
\begin{align} \label{eq:diff}
&\partial_\bsy^\bsnu\big[a\, \nabla u \cdot \nabla v + (\kappa\, u^2 + (f - \kappa)\,u)\,v\big]
\nonumber\\
&= \sum_{\bsk\le\bsnu} \binom{\bsnu}{\bsk}
  \Big[
  (\partial_\bsy^{\bsk} a)\,\nabla (\partial_\bsy^{\bsnu-\bsk} u)\cdot\nabla v 
  + (\partial_\bsy^\bsk \kappa)\,(\partial_\bsy^{\bsnu-\bsk} u^2)  \, v 
  + (\partial_\bsy^\bsk (f - \kappa))\,(\partial_\bsy^{\bsnu-\bsk} u) \, v
  \Big],
\end{align}%
where the mixed derivative of $u^2$ can be expanded further using the Leibniz product rule as
\begin{align*}
  \partial_\bsy^{\bsnu-\bsk} u^2
  &= \sum_{\bsm\le\bsnu-\bsk} 
  \binom{\bsnu-\bsk}{\bsm} (\partial_\bsy^\bsm u)\,(\partial_\bsy^{\bsnu-\bsk-\bsm} u), 
\end{align*}
for which the $\bsk=\bszero$ case can be written as
\begin{align*}
  \partial_\bsy^{\bsnu} u^2
  &= \sum_{\bsm\le\bsnu} 
  \binom{\bsnu}{\bsm} (\partial_\bsy^\bsm u)\,(\partial_\bsy^{\bsnu-\bsm} u) 
  = 2\,(\partial_\bsy^\bsnu u)\,u + \sum_{\bszero<\bsm <\bsnu} 
  \binom{\bsnu}{\bsm} (\partial_\bsy^\bsm u)\,(\partial_\bsy^{\bsnu-\bsm} u).
\end{align*}
Separating out the $\bsk=\bszero$ terms in \eqref{eq:diff} and applying \eqref{eq:diff-ak}, we obtain for all $v \in\calX$,
\begin{align} \label{eq:der}
  &\int_0^\tau\!\!
  \big\langle \tfrac{\partial}{\partial t} (\partial_\bsy^\bsnu u),v\big\rangle_{H^1(\Omega)^*,H^1(\Omega)}
  \,\rd t
  + \int_0^\tau\!\! \int_\Omega 
  \bigg[
  a\,\nabla (\partial_\bsy^\bsnu u)\cdot\nabla v
  + 2\,\kappa\,(\partial_\bsy^\bsnu u)\,u \, v
  + (f-\kappa)\,(\partial_\bsy^\bsnu u) \, v \\
  &\qquad
  + \sum_{j\ge 1} \nu_j\,\psi_j\,\nabla (\partial_\bsy^{\bsnu-\bse_j} u)\cdot\nabla v
  - \sum_{j\ge 1} \nu_j\,\xi_j\,(\partial_\bsy^{\bsnu-\bse_j} u) \, v 
  + \kappa\,\sum_{\bszero < \bsm <\bsnu} \binom{\bsnu}{\bsm}\, 
    (\partial_\bsy^\bsm u)\,(\partial_\bsy^{\bsnu-\bsm} u) \, v \nonumber\\
  &\qquad
  + \sum_{j\ge 1} \nu_j\,\xi_j \sum_{\bsm\le\bsnu-\bse_j} \binom{\bsnu-\bse_j}{\bsm}\, 
   (\partial_\bsy^\bsm u)\,(\partial_\bsy^{\bsnu-\bse_j-\bsm} u) \, v
     \bigg]\,\rd\bsx\,\rd t = 0. \nonumber
\end{align}

\textsc{Membership in $\calX$.} Taking $\tau = T$ in \eqref{eq:der} and moving all terms involving derivatives of order strictly less than $\bsnu$ to the right-hand side, we can see that
$\partial_\bsy^\bsnu u$ satisfies (the weak form of) a \emph{linear} parabolic PDE
\begin{align}
\label{eq:der_weak_eq}
&\int_I \big\langle \tfrac{\partial}{\partial t} (\partial_\bsy^\bsnu u),v\big\rangle_{H^1(\Omega)^*,H^1(\Omega)}
  \,\rd t
  + \int_I \int_\Omega \big[
  a\,\nabla (\partial_\bsy^\bsnu u)\cdot\nabla v + \widetilde{b}\, (\partial_\bsy^\bsnu u)\, v\big] \, \rd \bsx \,\rd t
  \\\nonumber
  &= \int_I \langle F_\bsnu, v\rangle_{H^1(\Omega)^*,H^1(\Omega)} \, \rd t, 
  \qquad \forall v \in \calX,
\end{align}
where $\widetilde{b} = \widetilde{b}^\bsy(\bsx, t) = 2\,\kappa\, u + f - \kappa \in L^2(I; L^\infty(\Omega))$ and $F_\bsnu = F_\bsnu(t)$ is the (time-dependent) functional defined by 
\begin{align}
\label{eq:F_nu}
& \langle F_\bsnu, v\rangle_{H^1(\Omega)^*,H^1(\Omega)} :=
- \int_\Omega \bigg[\sum_{j\ge 1} \nu_j\,\psi_j\,\nabla (\partial_\bsy^{\bsnu-\bse_j} u)\cdot\nabla v
  - \sum_{j\ge 1} \nu_j\,\xi_j\,(\partial_\bsy^{\bsnu-\bse_j} u) \, v 
  \\\nonumber
  + \kappa\,&\sum_{\bszero < \bsm <\bsnu} \binom{\bsnu}{\bsm}\, 
    (\partial_\bsy^\bsm u)\,(\partial_\bsy^{\bsnu-\bsm} u) \, v
  + \sum_{j\ge 1} \nu_j\,\xi_j \sum_{\bsm\le\bsnu-\bse_j} \binom{\bsnu-\bse_j}{\bsm}\, 
   (\partial_\bsy^\bsm u)\,(\partial_\bsy^{\bsnu-\bse_j-\bsm} u) \, v
     \bigg]\,\rd\bsx,
\end{align}
which will be shown to be well defined below. Since $u_0$ is independent of $\bsy$, differentiating the initial 
condition in \eqref{eq:weak1} gives the initial condition
$\int_\Omega (\partial_\bsy^\bsnu u(\cdot, 0))\, v \, \rd \bsx = 0$ for all $v \in L^2(\Omega)$.

We now prove by induction that $\partial_\bsy^\bsnu u\in \calX$ for all multiindices $\bsnu$. For the base case $\bsnu = \bszero$, we have $\partial_\bsy^\bszero u = u \in \calX$ by Theorem~\ref{thm:exist}. Then under the induction hypothesis that $\partial_\bsy ^\bsm u \in \calX$ for all $\bsm < \bsnu$, it follows that $F_\bsnu \in L^2(I; H^1(\Omega))^* = L^2(I; H^1(\Omega)^*)$, which we show below by considering functionals corresponding to each term on the right-hand side of \eqref{eq:F_nu} separately.

Let $v \in L^2(I; H^1(\Omega))$. For the first term in \eqref{eq:F_nu} involving gradients, 
using $\psi_j \in L^\infty(\Omega)$ and $\partial_\bsy^{\bsnu-\bse_j} u\in\calX$ then applying the Cauchy--Schwarz inequality twice gives
\begin{align*}
\int_I \int_\Omega \psi_j\,\nabla (\partial_\bsy^{\bsnu-\bse_j} u) \cdot \nabla v \, \rd \bsx \, \rd t 
&\leq \|\psi_j\|_{L^\infty(\Omega)} \int_I \|\nabla (\partial_\bsy^{\bsnu-\bse_j} u)\|_{L^2(\Omega)}\, \|\nabla v\|_{L^2(\Omega)} \, \rd t
\\
&\leq \|\psi_j\|_{L^\infty(\Omega)}\, \|\partial_\bsy^{\bsnu-\bse_j} u\|_{\calX} \, \| v\|_{L^2(I; H^1(\Omega))}.
\end{align*}
Hence, the functionals corresponding to the gradient terms in \eqref{eq:F_nu} are in $L^2(I; H^1(\Omega)^*)$.
Similarly, for the second term in \eqref{eq:F_nu} we obtain
\begin{align*}
\int_I \int_\Omega \xi_j\,(\partial_\bsy^{\bsnu-\bse_j} u) \, v \, \rd \bsx \, \rd t 
&\le \|\xi_j\|_{L^\infty(\Omega)}\, \|\partial_\bsy^{\bsnu-\bse_j} u\|_{\calX} \, \| v\|_{L^2(I; H^1(\Omega))}.
\end{align*}
For the third and fourth terms in \eqref{eq:F_nu}, we face a new challenge to manage the space-time integrals of three factors involving $u$.
If we blindly apply the Cauchy--Schwarz inequality multiple times to separate the three factors, then we would raise norms to higher powers than~$2$ and the subsequent induction would fail.

To handle integrals of the form $\int_I \int_{\Omega} f\,g\,h \,\rd\bsx\,\rd t$ where the integrand is a product of three factors and to keep the norms with exponent $2$, we apply the generalized H\"older inequality, first for $\bsx\in\Omega$ with $\frac{1}{2}+\frac{1}{4}+\frac{1}{4}=1$, and then for $t\in I$ with $\frac{1}{\infty}+\frac{1}{2}+\frac{1}{2}=1$, as follows:
\begin{align} \label{eq:fantastic}
  &\int_I \int_{\Omega} f\,g\,h \,\rd\bsx\,\rd t
  \le \int_I \|f\|_{L^2(\Omega)}\,\|g\|_{L^4(\Omega)}\,\|h\|_{L^4(\Omega)}\,\rd t \nonumber\\
  &\le \Big(\max_{t\in I} \|f\|_{L^2(\Omega)}\Big)\,
  \Big(\int_I \|g\|_{L^4(\Omega)}^2\,\rd t\Big)^{1/2}
  \Big(\int_I \|h\|_{L^4(\Omega)}^2\,\rd t\Big)^{1/2} \nonumber\\
  &\le \theta_1^2\,\Big(\max_{t\in I} \|f\|_{L^2(\Omega)}\Big)\,
  \Big(\int_I \|g\|_{H^1(\Omega)}^2\,\rd t\Big)^{1/2}
  \Big(\int_I \|h\|_{H^1(\Omega)}^2\,\rd t\Big)^{1/2} \nonumber\\
  &= \theta_1^2\,\|f\|_{C(I;L^2(\Omega))}\,\|g\|_{L^2(I;H^1(\Omega))}\,
  \|h\|_{L^2(I;H^1(\Omega))} 
  \le \theta_1^2\,\theta_2\,\|f\|_\calX\,\|g\|_\calX\,\|h\|_{L^2(I;H^1(\Omega))},
\end{align}
where we used the embeddings \eqref{eq:emb2} and \eqref{eq:emb3}. Thus we obtain for the third and fourth terms in \eqref{eq:F_nu},
\begin{align*}
    \int_I \int_\Omega \kappa\, (\partial_\bsy^\bsm u)\, (\partial_\bsy^{\bsnu - \bsm} u) v \, \rd \bsx \, \rd t
     & \leq \theta_1^2\,\theta_2\, \kappa_{\max}\,\|\partial_\bsy^\bsm u\|_\calX\,\|\partial_\bsy^{\bsnu - \bsm} u\|_\calX\,\|v\|_{L^2(I;H^1(\Omega))},
     \\
    \int_I \int_\Omega \xi_j\, (\partial_\bsy^\bsm u) (\partial_\bsy^{\bsnu - \bse_j - \bsm} u)\, v \, \rd \bsx \, \rd t
     & \leq \theta_1^2\,\theta_2\, \|\xi_j\|_{L^\infty(\Omega)}\,\|\partial_\bsy^\bsm u\|_\calX\,\|\partial_\bsy^{\bsnu - \bse_j - \bsm} u\|_\calX\,\|v\|_{L^2(I;H^1(\Omega))}.
\end{align*}
Seeing as though $F_\bsnu$ is a linear combination of functionals belonging to $L^2(I; H^1(\Omega)^*)$, we conclude that $F_\bsnu \in L^2(I; H^1(\Omega)^*)$.

Since the bilinear form corresponding to the second term in \eqref{eq:der_weak_eq} is elliptic and bounded and $F_\bsnu \in L^2(I; H^1(\Omega)^*)$, it then follows from standard results for linear parabolic PDEs (e.g., \cite{DL92,SchSt09}) that there exists a unique $\partial_\bsy^\bsnu u \in \calX$ satisfying \eqref{eq:der_weak_eq}.

\textsc{Gradient.}
Taking the mixed derivative $\partial_\bsy^\bsnu$ on both sides of the second condition in \eqref{eq:weak1} and then taking $v = \partial_\bsy^\bsnu u(\cdot,0)$, we conclude that $\|\partial_\bsy^\bsnu u(\cdot,0)\|_{L^2(\Omega)}^2 = 0$. Since $\partial_\bsy^\bsnu u \in \calX$, by \eqref{eq:integration_by_parts},
this yields
\begin{align*}
 \int_0^{\tau}\!\! 
 \big\langle \tfrac{\partial}{\partial t} (\partial_\bsy^\bsnu u),
 \partial_\bsy^\bsnu u\big\rangle_{H^1(\Omega)^*,H^1(\Omega)}
 \,\rd t 
 = \int_0^{\tau}\! 
 \frac{1}{2} \frac{\rd}{\rd t} \|\partial_\bsy^\bsnu u(\cdot,t)\|_{L^2(\Omega)}^2\rd t 
 = \frac{\|\partial_\bsy^\bsnu u(\cdot,\tau)\|_{L^2(\Omega)}^2}{2}
 \ge 0.
\end{align*}
Taking $v = \partial_\bsy^\bsnu u$ in \eqref{eq:der} then gives
\begin{align} \label{eq:equal} 
  \int_0^\tau\! \frac{1}{2} \frac{\rd}{\rd t} \|\partial_\bsy^\bsnu u(\cdot,t)\|_{L^2(\Omega)}^2 \rd t
  + \!\int_0^\tau\!\! \int_\Omega 
  \Big[ a\,|\nabla (\partial_\bsy^\bsnu u)|^2
  + 2\,\kappa\,(\partial_\bsy^\bsnu u)^2\,u 
  + (f-\kappa)\,(\partial_\bsy^\bsnu u)^2 + \eta \Big] \rd\bsx\,\rd t 
  = 0, 
\end{align}
where
\begin{align} \label{eq:eta}
  \eta(\bsx,t) 
  &:= \sum_{j\ge 1} \nu_j\,\psi_j\,
  \nabla (\partial_\bsy^{\bsnu-\bse_j} u)\cdot\nabla (\partial_\bsy^\bsnu u)
  - \sum_{j\ge 1} \nu_j\,\xi_j\,(\partial_\bsy^{\bsnu-\bse_j} u) \, (\partial_\bsy^\bsnu u) \\
  &\quad + \kappa \sum_{\bszero < \bsm <\bsnu} \binom{\bsnu}{\bsm}\,
  (\partial_\bsy^\bsm u)\,(\partial_\bsy^{\bsnu-\bsm} u)\,(\partial_\bsy^\bsnu u) \nonumber\\
  &\quad + \sum_{j\ge 1} \nu_j\,\xi_j \sum_{\bsm\le\bsnu-\bse_j} \binom{\bsnu-\bse_j}{\bsm}\, 
  (\partial_\bsy^\bsm u)\,(\partial_\bsy^{\bsnu-\bse_j-\bsm} u) \, (\partial_\bsy^\bsnu u). \nonumber
\end{align}
Note that $\eta$ can potentially give both positive and negative values.%

Until this point we still have equality in \eqref{eq:equal}--\eqref{eq:eta}. For the special cases $|\bsnu|=1$, the sum over $\bszero < \bsm <\bsnu$ is empty and takes the value $0$, with no meaning attached to the summand.

Recall that $0\le u\le 1$. Now we drop the nonnegative terms $\kappa\,(\partial_\bsy^\bsnu u)^2\,u$ and $f\,(\partial_\bsy^\bsnu u)^2$ in \eqref{eq:equal} and rearrange to obtain the inequality
\begin{align*}  
  \int_0^\tau\! \frac{1}{2} \frac{\rd}{\rd t} \|\partial_\bsy^\bsnu u(\cdot,t)\|_{L^2(\Omega)}^2\,\rd t 
  + \int_0^\tau\!\! \int_\Omega \Big[a\,|\nabla (\partial_\bsy^\bsnu u)|^2 + \eta\Big]\,\rd\bsx \,\rd t 
  \le \kmax \int_0^\tau\! \|\partial_\bsy^\bsnu u\|_{L^2(\Omega)}^2 \,\rd t.
\end{align*}
If we were to return to the derivation of \eqref{eq:weak1}, but multiply the PDE by $e^{-2\,\kmax\,t}$ before integrating over $t\in [0,\tau]$, then the preceding argument would lead to the above inequality but with all integrands multiplied by $e^{-2\,\kmax\,t}$, and by the product rule it is equivalent to
\begin{align} \label{eq:wow}
  \int_0^\tau\! \frac{1}{2} \frac{\rd}{\rd t} 
  \Big(\|\partial_\bsy^\bsnu u(\cdot,t)\|_{L^2(\Omega)}^2 e^{-2\,\kmax\, t}\Big)\,\rd t 
  + \int_0^\tau\!\! \int_\Omega \Big[a\,|\nabla (\partial_\bsy^\bsnu u)|^2 + \eta\Big] e^{-2\,\kmax\, t}
  \,\rd\bsx\,\rd t
  \le 0. 
\end{align}

Using $\|\partial_\bsy^\bsnu u(\cdot,0)\|_{L^2(\Omega)}^2 = 0$ and \eqref{eq:integration_by_parts}, we drop the nonnegative term $a\,|\nabla (\partial_\bsy^\bsnu u)|^2$ in \eqref{eq:wow} to obtain 
\begin{align*} 
    \frac{\|\partial_\bsy^\bsnu u(\cdot,\tau)\|_{L^2(\Omega)}^2 e^{-2\,\kmax\,\tau} }{2} = 
  \int_0^\tau\! \frac{1}{2} \frac{\rd}{\rd t} 
  \Big(\|\partial_\bsy^\bsnu u(\cdot,t)\|_{L^2(\Omega)}^2 e^{-2\,\kmax\, t}\Big)\,\rd t 
  \le \int_0^\tau\!\! \int_\Omega |\eta| \,e^{-2\,\kmax\,t}\,\rd\bsx\,\rd t.
\end{align*}
We multiply both sides of the inequality by $2\,e^{2\,\kmax\,\tau}$ and integrate over $\tau\in I$ to obtain
\begin{align} \label{eq:apple1}
 \|\partial_\bsy^\bsnu u\|_{L^2(I;L^2(\Omega))}^2
 &\le 2\int_0^T\!\! \int_0^\tau\!\!\int_\Omega |\eta|\,e^{2\,\kmax\,(\tau-t)} \,\rd\bsx\,\rd t\,\rd\tau 
 = \int_I\int_\Omega |\eta|\,\frac{e^{2\,\kmax\, (T-t)}-1}{\kmax}\,\rd\bsx\,\rd t.
\end{align}
Next we drop the nonnegative first term in \eqref{eq:wow} and take $\tau = T$ to obtain
\begin{align} \label{eq:apple2}
 \amin\,e^{-2\,\kmax\,T}\,\|\nabla (\partial_\bsy^\bsnu u)\|_{L^2(I;L^2(\Omega))}^2
 \le \int_I \int_\Omega |\eta| \,e^{-2\,\kmax\,t}\,\rd\bsx\,\rd t.
\end{align}
Combining \eqref{eq:apple1} and \eqref{eq:apple2}, we conclude that
\begin{align} \label{eq:apple3}
 &\|\partial_\bsy^\bsnu u\|_{L^2(I;H^1(\Omega))}^2
 = \|\partial_\bsy^\bsnu u\|_{L^2(I;L^2(\Omega))}^2
 + \|\nabla (\partial_\bsy^\bsnu u)\|_{L^2(I;L^2(\Omega))}^2 \nonumber\\
 &\le \int_I \int_\Omega |\eta|\, 
 \Big(\frac{e^{2\,\kmax\, (T-t)}-1}{\kmax}  + \frac{e^{2\,\kmax\, (T-t)}}{\amin}\Big)
 \,\rd\bsx\,\rd t 
 \le \frac{2\,e^{2\,\kmax\,T}}{\min(\amin,\kmax)}
 \int_I \int_\Omega |\eta| \,\rd\bsx\,\rd t.
\end{align}%

To bound the space-time integral of $|\eta|$ with $\eta$ defined in \eqref{eq:eta}, we pull out $\beta_j := \max(\|\psi_j\|_\infty,\|\xi_j\|_\infty)$ and then apply the Cauchy--Schwarz inequality to the first two terms, as is typically done.
For the last two terms in \eqref{eq:eta} we again use \eqref{eq:fantastic} to manage the space-time integrals of three factors.

We obtain from \eqref{eq:eta}, \eqref{eq:apple3}, and \eqref{eq:fantastic},
\begin{align*} 
 &\frac{\min(\amin,\kmax)}{2\,e^{2\,\kmax\,T}}\,
 \|\partial_\bsy^\bsnu u\|_{L^2(I;H^1(\Omega))}^2 
 \le \int_I \int_\Omega |\eta| \,\rd\bsx\,\rd t \\
 &\le \sum_{j\ge 1} \nu_j\, \beta_j\,
  \|\nabla(\partial_\bsy^{\bsnu-\bse_j} u)\|_{L^2(I;L^2(\Omega))}\,
  \|\nabla(\partial_\bsy^\bsnu u)\|_{L^2(I;L^2(\Omega))} \nonumber\\
  &\quad + \sum_{j\ge 1} \nu_j\, \beta_j\,
  \|\partial_\bsy^{\bsnu-\bse_j} u\|_{L^2(I;L^2(\Omega))}\,
  \|\partial_\bsy^\bsnu u\|_{L^2(I;L^2(\Omega))} \nonumber\\
  &\quad + \theta_1^2\,\theta_2\,\kmax \sum_{\bszero\ne\bsm <\bsnu} \binom{\bsnu}{\bsm}\,
  \|\partial_\bsy^\bsm u\|_\calX\,
  \|\partial_\bsy^{\bsnu-\bsm} u\|_\calX\,
  \|\partial_\bsy^\bsnu u\|_{L^2(I;H^1(\Omega))} \nonumber\\
  &\quad + \theta_1^2\,\theta_2 \sum_{j\ge 1} \nu_j\,\beta_j\,
  \sum_{\bsm\le\bsnu-\bse_j} 
  \binom{\bsnu-\bse_j}{\bsm}
  \|\partial_\bsy^\bsm u\|_\calX\,
  \|\partial_\bsy^{\bsnu-\bse_j-\bsm} u\|_\calX\,
  \|\partial_\bsy^\bsnu u\|_{L^2(I;H^1(\Omega))}.
\end{align*} 
The last factor in every term on the right-hand side can be upper bounded by $\|\partial_\bsy^\bsnu u\|_{L^2(I;H^1(\Omega))}$. Canceling out this common factor from both sides and grouping the first two terms on the right-hand side, we arrive at
\begin{align} \label{eq:part1}
  &\frac{\min(\amin,\kmax)}{2\,e^{2\,\kmax\,T}}\,
 \|\partial_\bsy^\bsnu u\|_{L^2(I;H^1(\Omega))} \nonumber\\
  &\le 2\sum_{j\ge 1} \nu_j\,\beta_j\,
  \|\partial_\bsy^{\bsnu-\bse_j} u\|_\calX 
  + \theta_1^2\,\theta_2\,\kmax
  \sum_{\bszero\ne\bsm <\bsnu} \binom{\bsnu}{\bsm}\, 
  \|\partial_\bsy^\bsm u\|_\calX\, 
  \|\partial_\bsy^{\bsnu-\bsm} u\|_\calX \nonumber\\
  &\qquad + \theta_1^2\,\theta_2 \sum_{j\ge 1} \nu_j\,\beta_j
  \sum_{\bsm\le\bsnu-\bse_j} \binom{\bsnu-\bse_j}{\bsm}\,
  \|\partial_\bsy^\bsm u\|_\calX\, \|\partial_\bsy^{\bsnu-\bse_j-\bsm} u\|_\calX
  \,=:\, \Theta_{\bsnu}.
\end{align} 

\textsc{Time derivative.}
We return to \eqref{eq:der} but remove the integration over $t$ and move all terms except the duality pair to the right-hand side. Then we apply the triangle, Cauchy--Schwarz and H\"older inequalities, together with $|2\,\kappa\,(\partial_\bsy^\bsnu u)\,u \, v + (f-\kappa)\,(\partial_\bsy^\bsnu u)\,v |\le (\kappa + f)\,|\partial_\bsy^\bsnu u|\,|v|$ since $0\le u\le 1$. We obtain for a.e.~$t\in I$,
\begin{align*} 
  &\big|
  \big\langle\tfrac{\partial}{\partial t} (\partial_\bsy^\bsnu u),v\big\rangle_{H^1(\Omega)^*,H^1(\Omega)}
  \big| \\
  &\le \amax\, \|\nabla (\partial_\bsy^\bsnu u)(\cdot,t)\|_{L^2(\Omega)}\,\|\nabla v\|_{L^2(\Omega)}  
  + (\kmax + \fmax)\, \|\partial_\bsy^\bsnu u(\cdot,t)\|_{L^2(\Omega)}\, \|v\|_{L^2(\Omega)} \\
  &\quad + \sum_{j\ge 1} \nu_j\, \beta_j\,
  \|\nabla (\partial_\bsy^{\bsnu-\bse_j} u)(\cdot,t)\|_{L^2(\Omega)}\,
  \|\nabla v\|_{L^2(\Omega)} 
  + \sum_{j\ge 1} \nu_j\, \beta_j\,
  \|\partial_\bsy^{\bsnu-\bse_j} u(\cdot,t)\|_{L^2(\Omega)}\,
  \|v\|_{L^2(\Omega)} \\
  &\quad + \kmax \sum_{\bszero < \bsm <\bsnu} \binom{\bsnu}{\bsm}\,
  \|\partial_\bsy^\bsm u(\cdot,t)\|_{L^2(\Omega)}\,
  \|\partial_\bsy^{\bsnu-\bsm} u(\cdot,t)\|_{L^4(\Omega)}\, \|v\|_{L^4(\Omega)}\\
  &\quad + 
  \sum_{j\ge 1} \nu_j\,\beta_j \sum_{\bsm\le\bsnu-\bse_j} \binom{\bsnu-\bse_j}{\bsm}\,
  \|\partial_\bsy^\bsm u(\cdot,t)\|_{L^2(\Omega)}\,
  \|\partial_\bsy^{\bsnu-\bse_j-\bsm} u(\cdot,t)\|_{L^4(\Omega)}\, \|v\|_{L^4(\Omega)}.
\end{align*} 
Using $\|\nabla v\|_{L^2(\Omega)}\le \|v\|_{H^1(\Omega)}$ and $\|v\|_{L^4(\Omega)}\le \theta_1\,\|v\|_{H^1(\Omega)}$, the last factor of every term on the right-hand side can be bounded in terms of $\|v\|_{H^1(\Omega)}$. Therefore we conclude that
\begin{align} \label{eq:tosquare} 
  \|\tfrac{\partial}{\partial t} (\partial_\bsy^\bsnu u)(\cdot,t)\|_{H^1(\Omega)^*} 
  &\le (\amax+\kmax+\fmax)\, \|\partial_\bsy^\bsnu u(\cdot,t)\|_{H^1(\Omega)} 
  + 2\sum_{j\ge 1} \nu_j\, \beta_j\,
  \|\partial_\bsy^{\bsnu-\bse_j} u(\cdot,t)\|_{H^1(\Omega)} \nonumber\\
  &\quad + \theta_1^2\,\theta_2\,\kmax \sum_{\bszero\ne\bsm <\bsnu} \binom{\bsnu}{\bsm}\,
  \|\partial_\bsy^\bsm u\|_\calX\,
  \|\partial_\bsy^{\bsnu-\bsm} u(\cdot,t)\|_{H^1(\Omega)} \nonumber\\
  &\quad + \theta_1^2\,\theta_2\,
  \sum_{j\ge 1} \nu_j\,\beta_j \sum_{\bsm\le\bsnu-\bse_j} \binom{\bsnu-\bse_j}{\bsm}\,
  \|\partial_\bsy^\bsm u\|_\calX\,
  \|\partial_\bsy^{\bsnu-\bse_j-\bsm} u(\cdot,t)\|_{H^1(\Omega)}, 
\end{align} 
where we also used $\|\partial_\bsy^\bsm u(\cdot,t)\|_{L^2(\Omega)}\le \|\partial_\bsy^\bsm u\|_{C(I;L^2(\Omega))} \le \theta_2\,\|\partial_\bsy^\bsm u\|_\calX$. This step to remove the dependence of one factor on $t$ is crucial before we integrate over $t\in I$.

Next we take the $L^2(I)$ norm of \eqref{eq:tosquare} and apply the triangle inequality to obtain
\begin{align} \label{eq:part2}
  \|\tfrac{\partial}{\partial t} (\partial_\bsy^\bsnu u)\|_{L^2(I;H^1(\Omega)^*)} 
  &\le (\amax+\kmax+\fmax)\, \|\partial_\bsy^\bsnu u\|_{L^2(I;H^1(\Omega))} 
  + \Theta_\bsnu, 
\end{align} 
where $\Theta_\bsnu$ is precisely the right-hand side of \eqref{eq:part1}.

\textsc{Recursion.}
Using \eqref{eq:part1} and \eqref{eq:part2}, we obtain
\begin{align*} 
  \|\partial_\bsy^\bsnu u\|_\calX
  &= \sqrt{ 
  \|\partial_\bsy^\bsnu u\|_{L^2(I;H^1(\Omega))}^2
  +\|\tfrac{\partial}{\partial t} (\partial_\bsy^\bsnu u)\|_{L^2(I;H^1(\Omega)^*)}^2} \\
  &\le \|\partial_\bsy^\bsnu u\|_{L^2(I;H^1(\Omega))} 
  + \|\tfrac{\partial}{\partial t} (\partial_\bsy^\bsnu u)\|_{L^2(I;H^1(\Omega)^*)} \\
  &\le (1 + \amax+\kmax+\fmax)\, \|\partial_\bsy^\bsnu u\|_{L^2(I;H^1(\Omega))} 
  + \Theta_\bsnu \\
  &\le (1 + \amax+\kmax+\fmax)\, \frac{2\,e^{2\kmax\,T}\,\Theta_\bsnu}{\min(\amin,\kmax)}
  + \Theta_\bsnu,
\end{align*}
which yields the recursion
\begin{align} \label{eq:recur}
  \|\partial_\bsy^\bsnu u\|_\calX
  &\le \Phi\, \bigg[ 2\sum_{j\ge 1} \nu_j\,\beta_j\,\|\partial_\bsy^{\bsnu-\bse_j} u\|_\calX 
  + \theta_1^2\,\theta_2\,\kmax \sum_{\bszero\ne\bsm <\bsnu} \binom{\bsnu}{\bsm}\, 
  \|\partial_\bsy^\bsm u\|_\calX\, 
  \|\partial_\bsy^{\bsnu-\bsm} u\|_\calX \nonumber\\
  &\qquad\quad + \theta_1^2\,\theta_2 \sum_{j\ge 1} \nu_j\,\beta_j
  \sum_{\bsm\le\bsnu-\bse_j} \binom{\bsnu-\bse_j}{\bsm}\,
  \|\partial_\bsy^\bsm u\|_\calX\, 
  \|\partial_\bsy^{\bsnu-\bse_j-\bsm} u\|_\calX \bigg],
\end{align}
where we abbreviated
\begin{align} \label{eq:Phi}
  \Phi := \frac{2\,e^{2\,\kmax\,T} (1+\amax+\kmax+\fmax)}{\min(\amin,\kmax)} + 1.
\end{align}

\textsc{Induction.}
We will prove by induction on $|\bsnu|$ that
\begin{align} \label{eq:hyp}
  \|\partial_\bsy^\bsnu u\|_\calX
  \le D\, \bsbeta^{\bsnu}\rho^{|\bsnu|-1}\, \bhalf_{|\bsnu|}
  \qquad\mbox{for all } \bsnu\ne\bszero,
\end{align} 
with constants $D$ and $\rho$ to be specified below.

Our base step is $|\bsnu|=1$ (rather than $\bsnu=\bszero$, as in typical proofs of parametric regularity) because our recursion \eqref{eq:recur} involves a sum over $\bszero < \bsm < \bsnu$ which is empty if $|\bsnu|=1$. Recall from Theorem~\ref{thm:priori} that we have the a priori bound $\|u\|_\calX \le C\,\|u_0\|_{L^2(\Omega)}$, where $\|u_0\|_{L^2(\Omega)}^2\le |\Omega|$ since $0\le u_0(\bsx)\le 1$.
For $\bsnu=\bse_k$ for some $k\ge 1$, the recursion \eqref{eq:recur} gives
\begin{align*} 
  \|\partial_\bsy^{\bse_k} u\|_\calX
  &\le \Phi\, \big( 2\,\beta_k\,\|u\|_\calX 
  + \theta_1^2\,\theta_2\, \beta_k\, \|u\|_\calX^2\big) \\
  &\le \Phi\, \big( 2\,\beta_k\,C\,\|u_0\|_{L^2(\Omega)}
  + \theta_1^2\,\theta_2\,\beta_k\, C^2\,\|u_0\|_{L^2(\Omega)}\,|\Omega|^{1/2}\big)
  \le \frac{D\,\beta_k}{2},
\end{align*}
as required in \eqref{eq:hyp}, if we define 
\begin{align} \label{eq:DE}
  D := E\, \|u_0\|_{L^2(\Omega)}, 
  \quad\mbox{with}\quad
  E := 2\,\Phi\,C\,\big( 2 +  C\,\theta_1^2\,\theta_2\,|\Omega|^{1/2}\big).
\end{align}

For the induction step, let $|\bsnu|\ge 2$ and suppose \eqref{eq:hyp} is true for all multiindices of order less than $|\bsnu|$. 
We have from \eqref{eq:recur}
\begin{align*} 
 &(1/\Phi)\,\|\partial_\bsy^\bsnu u\|_\calX \\
 &\le 2 \sum_{j\ge 1} \nu_j\, \beta_j\,
  D\,\bsbeta^{\bsnu-\bse_j}\rho^{|\bsnu-\bse_j|-1}\,
  \bhalf_{|\bsnu-\bse_j|} \\
  &\quad + \theta_1^2\,\theta_2\,\kmax
  \sum_{\bszero\ne\bsm <\bsnu} \binom{\bsnu}{\bsm} 
  D\,\bsbeta^{\bsm}\rho^{|\bsm|-1}\,\bhalf_{|\bsm|}\,
  D\,\bsbeta^{\bsnu-\bsm}\rho^{|\bsnu-\bsm|-1}\, \bhalf_{|\bsnu-\bsm|} \\
  &\quad 
  + \theta_1^2\,\theta_2
  \sum_{j\ge 1} \nu_j\,\beta_j \!\!\!\!
  \sum_{\bszero\ne\bsm < \bsnu-\bse_j} \!\!\! \binom{\bsnu-\bse_j}{\bsm}
  D\,\bsbeta^{\bsm}\rho^{|\bsm|-1}\,\bhalf_{|\bsm|}\,
  D\,\bsbeta^{\bsnu-\bse_j-\bsm}\rho^{|\bsnu-\bse_j-\bsm|-1}\,
  \bhalf_{|\bsnu-\bse_j-\bsm|} \\
  &\quad + 2\,\theta_1^2\,\theta_2 \sum_{j\ge 1} \nu_j\,\beta_j\,
  D\,D\,\bsbeta^{\bsnu-\bse_j}\rho^{|\bsnu-\bse_j|-1}\,\bhalf_{|\bsnu-\bse_j|}
  \\
  &\le 2\,D\,\bsbeta^{\bsnu}\rho^{|\bsnu|-2}\, |\bsnu| \,\bhalf_{|\bsnu|-1}
  + D^2\theta_1^2\,\theta_2\,\kmax\,\bsbeta^{\bsnu}\rho^{|\bsnu|-2}
  \sum_{\bszero\ne\bsm <\bsnu} \binom{\bsnu}{\bsm}\, \bhalf_{|\bsm|}\, \bhalf_{|\bsnu-\bsm|} \\
  &\quad + D^2\theta_1^2\,\theta_2\,\bsbeta^{\bsnu}\rho^{|\bsnu|-3}\,
  \sum_{j\ge 1} \nu_j\sum_{\bszero\ne\bsm<\bsnu-\bse_j} \binom{\bsnu-\bse_j}{\bsm}\,
  \bhalf_{|\bsm|}\, \bhalf_{|\bsnu-\bse_j-\bsm|} \\
  &\quad + 2\,D^2\,\theta_1^2\,\theta_2\,
   \bsbeta^{\bsnu}\rho^{|\bsnu|-2}\, |\bsnu| \,\bhalf_{|\bsnu|-1},
\end{align*}
where we separated out the $\bsm=\bszero$ and $\bsm=\bsnu-\bse_j$ terms in the sum over $\bsm\le\bsnu-\bse_j$ and used $\|u\|_\calX \le C\,\|u_0\|_{L^2(\Omega)} \le D$.
Using \eqref{eq:falling}, this leads to
\begin{align*} 
  (1/\Phi)\,\|\partial_\bsy^\bsnu u\|_\calX 
  &\le 2\,D\,\bsbeta^{\bsnu}\rho^{|\bsnu|-2}\, |\bsnu|\,\bhalf_{|\bsnu|-1}
  + D^2\theta_1^2\,\theta_2\,\kmax\,\bsbeta^{\bsnu}\rho^{|\bsnu|-2}
  \cdot 2\,\bhalf_{|\bsnu|} \\
  &\quad + D^2\theta_1^2\,\theta_2\,\bsbeta^{\bsnu}\rho^{|\bsnu|-3}\,
  |\bsnu| \cdot 2\,\bhalf_{|\bsnu|-1}
  + 2\,D^2\,\theta_1^2\,\theta_2\,
   \bsbeta^{\bsnu}\rho^{|\bsnu|-2}\, |\bsnu| \,\bhalf_{|\bsnu|-1},\\
  &\le D\,\bsbeta^{\bsnu}\rho^{|\bsnu|-2}\,\bhalf_{|\bsnu|}\,
  \big(8 + 2E\,\theta_1^2\,\theta_2\,\kmax\,|\Omega|^{1/2}
  + 16E\,\theta_1^2\,\theta_2\,|\Omega|^{1/2}\big),
\end{align*}
where we assumed $\rho\ge 1$ so that $\rho^{|\bsnu|-3}\le \rho^{|\bsnu|-2}$, and used $D^2\le DE\,|\Omega|^{1/2}$ (see \eqref{eq:DE}), together with 
\begin{align*}
  n\,\bhalf_{n-1} 
  = \frac{n}{|\tfrac{1}{2}-n+1|}\,\bhalf_{n}
  \le \frac{n}{\frac{1}{4}n}\,\bhalf_{n}
  = 4\,\bhalf_{n}
  \quad\mbox{for all } n\ge 2.
\end{align*}
Defining
\begin{align} \label{eq:rho}
 \rho := 2\,\Phi\,
 \big(4 + E\,\theta_1^2\,\theta_2\,(\kmax + 8)\,|\Omega|^{1/2}\big) \ge 1,
\end{align}
we arrive at the required bound \eqref{eq:hyp}. This completes the induction proof.

We have written the hypothesis \eqref{eq:hyp} in a convenient form for induction, where the factor $\rho$ served as a scaling parameter that can be adjusted to suit our needs, as we defined in \eqref{eq:rho}. 
Since $D = E\, \|u_0\|_{L^2(\Omega)}$ (see \eqref{eq:DE}), this factor $\rho$ can be seen as a rescaling of the sequence $\beta_j$ as follows:
\begin{align*} 
  \|\partial_\bsy^\bsnu u\|_\calX
  \le D\,\bsbeta^{\bsnu}\rho^{|\bsnu|-1}\, \bhalf_{|\bsnu|}
  = \tfrac{D}{\rho}\, (\rho\,\bsbeta)^{\bsnu}\, \bhalf_{|\bsnu|}
  = \tfrac{E}{\rho}\, (\rho\,\bsbeta)^{\bsnu}\, \bhalf_{|\bsnu|}\,\|u_0\|_{L^2(\Omega)} 
  \quad\mbox{for all } \bsnu\ne\bszero.
\end{align*} 
The result can be extended to cover also the $\bsnu=\bszero$ case as follows
 \begin{align*} 
  \|\partial_\bsy^\bsnu u\|_\calX
  \le \max(C,\tfrac{E}{\rho})\, 
  (\rho\,\bsbeta)^{\bsnu}\,\bhalf_{|\bsnu|}\,\|u_0\|_{L^2(\Omega)} 
  \quad\mbox{for all } \bsnu.
\end{align*} 
This completes the proof.
\end{proof}

\section{QMC error analysis for uniform random fields} \label{sec:error}

As noted in the introduction, QMC methods have demonstrated impressive performance in several UQ applications involving PDEs with random coefficients or domains, including linear elliptic \glspl*{pde} \cite{dick2014higher,gantner2018,graham2015quasi,graham2018circulant,guth2021quasi,hakula2024,kuo2012quasi},
elliptic eigenvalue problems \cite{gilbert2019evp,nguyen2024evp}
and linear parabolic equations with control \cite{GKKSS24}.
For a more pedagogical survey of some of the earlier works from the QMC perspective, we refer the interested reader to~\cite{kuo2016application}.

We now fix $\psi_j \equiv 0$ for $j$ even and $\xi_j \equiv 0$ for $j$ odd to make the two fields independent.
We truncate the random fields to $s$ terms and denote the corresponding PDE solution by $u_s^{\bsy}$ with $\bsy = (y_1,\ldots,y_s, 0,\ldots)$.
We are interested in the expected value of the solution $u_s$ with respect to the random parameter $\bsy\in U_s := [-\frac{1}{2},\frac{1}{2}]^s$,
\begin{equation*} 
    \bbE[u_s] := \int_{U_s} u_s^\bsy\,\rd\bsy.
\end{equation*}
We approximate this integral using a randomly shifted lattice rule with generating vector $\bsz\in \bbZ^s$ and uniformly generated random shift $\bsDelta\in [0,1]^s$,
\begin{equation} \label{eq:qmc-estimate}
    Q_N(u_s) := \frac{1}{N} \sum_{i=1}^N u^{\bsy^{(i)}_{\bsDelta}},
    \quad\mbox{with}\quad
    \bsy^{(i)}_{\bsDelta} 
    := (\tfrac{i}{N}\,\bsz + \bsDelta) \bmod 1 - \bshalf,
\end{equation}
where the subtraction by $\bshalf=(\tfrac{1}{2},\ldots,\tfrac{1}{2})$ translates points from the standard unit cube $[0,1]^s$ to~$U_s$.
Alternatively, we may wish to compute the expectation of a QoI given by a bounded linear functional $G$ of the PDE solution $u_s^\bsy$, in which case we approximate $\bbE[G(u_s)] \approx Q_N(G(u_s^\bsy))$, with the QMC rule applied analogously as in \eqref{eq:qmc-estimate}. A randomly shifted lattice rule is unbiased, i.e., $\bbE_\bsDelta [Q_N(G(u_s))] = \bbE [G(u_s)]$.
For now we work with the exact PDE solution (i.e., no finite element method or time discretization).

Using the parametric regularity bound established in Theorem~\ref{thm:reg}, we can apply two general theorems from \cite{GKKSS24} for Banach-space-valued functions depending on affine uniform parameters $\bsy$, to obtain a dimension truncation error bound and a QMC error bound. We analyze the QMC error in the first-order weighted unanchored Sobolev space with weight parameters $\bsgamma := (\gamma_\setu)_{\setu\subset\bbN}$, see e.g., \cite{kuo2012quasi}.

\begin{theorem}[dimension truncation error]
\label{thm:trunc-error}
Suppose that $\|\psi_1\|_\infty\ge \|\psi_3\|_\infty\ge \|\psi_5\|_\infty\ge \cdots$ and $\|\xi_2\|_\infty\ge \|\xi_4\|_\infty\ge \|\xi_6\|_\infty\ge \cdots$, and 
$\sum_{j\ge 1} \beta_j^p = \sum_{j\ge 1} \max(\|\psi_j\|_\infty, \|\xi_j\|_\infty)^p < \infty$ for some $p\in (0,1)$. Then for all $s\in\bbN$ the dimension truncation error satisfies
\begin{align*}
 \big\|\,\bbE[u - u_s]\,\big\|_\calX 
 = \Big\|\int_U (u^\bsy - u^{\bsy_s})\,\rd\bsy \Big\|_\calX
 = \widetilde{C}\, \|u_0\|_{L^2(\Omega)}\, C^*\,s^{-2/p+1},
\end{align*}
where $\widetilde{C}$ is the constant from Theorem~\ref{thm:reg} and $C^*>0$ is independent of $s$.
\end{theorem}

\begin{proof}
Noting that $\bhalf_{|\bsnu|}\le |\bsnu|!$\,, we apply \cite[Theorem~6.2]{GKKSS24} by setting there $C_0 = \widetilde{C}\, \|u_0\|_{L^2(\Omega)}$, $r_1=0$, $r_2 = \rho$, and $\bsrho = \bsbeta$. Our interleaved sequence $\bsbeta$ here is not strictly ordered as required in \cite[Theorem~6.2]{GKKSS24}, but the same result holds with a simple modification of the proof by splitting tail sums into separate and ordered sums of odd and even terms.
\end{proof}

\begin{theorem}[QMC error bound]
\label{thm:qmc-error}
We have the root-mean-square error bound with respect to the random shift
\begin{align*}
 \sqrt{\bbE_\bsDelta\, \big[\, \big\| \bbE[u_s] - Q_N(u_s)\big\|_\calX^2 \big]}
 \le C_{s,\bsgamma,\lambda}\,[\phi_{\rm tot}(N)]^{-1/(2\lambda)}
 \quad\mbox{for all } \lambda\in (\tfrac{1}{2},1],
\end{align*}
or for any bounded linear functional $G : \calX\to \bbR$,
\begin{align*}
 \sqrt{\bbE_\bsDelta\, \big[\, \big| \bbE[G(u_s)] - Q_N(G(u_s)) \big|^2 \big]}
 \le \|G\|_{\calX^*}\, C_{s,\bsgamma,\lambda}\,[\phi_{\rm tot}(N)]^{-1/(2\lambda)}
 \quad\mbox{for all } \lambda\in (\tfrac{1}{2},1],
\end{align*}
where $\phi_{\rm tot}(N) := |\{1\le z\le N-1: \gcd(z,N)=1\}|$ is the Euler totient function, with 
$1/\phi_{\rm tot}(N) \le 2/N$ when $N$ is a prime power, and
\begin{align*}
  C_{s,\bsgamma,\lambda}
  := \widetilde{C}\,
  \bigg(\sum_{\emptyset\ne\setu\subseteq\{1:s\}} \gamma_\setu^\lambda 
  \Big(\frac{2\,\zeta(2\lambda)}{(2\pi^2)^\lambda}\Big)^{|\setu|}\bigg)^{1/(2\lambda)}
  \bigg(\sum_{\setu\subseteq\{1:s\}} 
  \frac{\bhalf_{|\setu|}^2\,\prod_{j\in\setu} (\rho\,\beta_j)^2}{\gamma_\setu} \bigg)^{1/2},
\end{align*}
where $\widetilde{C}$ is defined in Theorem~\ref{thm:reg} and $\rho$ is defined in \eqref{eq:rho}.

For some $p\in (0,1)$ suppose that $\sum_{j\ge 1} \beta_j^p = \sum_{j\ge 1} \max(\|\psi_j\|_\infty, \|\xi_j\|_\infty)^p < \infty$. Then the best convergence rate $\calO(N^{-\min(1-\delta,\,1/p-1/2)})$ is obtained by the following choices
\begin{align*}
  \lambda :=
  \begin{cases}
  \frac{1}{2-2\delta} \mbox{ for all } \delta\in (0,1) & \mbox{if } p\in [0,\frac{2}{3}], \\
  \frac{p}{2-p} & \mbox{if } p\in (\frac{2}{3},1),
  \end{cases}
  \quad
  \gamma_\setu := \bigg( \bhalf_{|\setu|}\,\prod_{j\in\setu} \frac{\rho\,\beta_j}{\sqrt{2\zeta(2\lambda)/(2\pi^2)^\lambda}} \bigg)^{2/(1+\lambda)}.
\end{align*}
and the constant $C_{s,\bsgamma,\lambda}$ is bounded independently of $s$.
\end{theorem}

\begin{proof}
First we adapt \cite[Theorem~6.6]{GKKSS24}, with $|\setu|!$ replaced by $\bhalf_{|\setu|}$ and setting $C_0 = \widetilde{C}\, \|u_0\|_{L^2(\Omega)}$, $r_1=0$, $r_2 = \rho$, $\rho_j = \beta_j$, to obtain the root-mean-square error bounds. Then we follow \cite[Theorem~6.8]{GKKSS24} to choose the parameter $\lambda$ and the function space weights $\gamma_\setu$. 
\end{proof}

\section{Numerical experiments}
\label{sec:numerics}

\subsection{Experiments with uniform random fields}
\label{sec:uniform-numerics}
We consider the spatial domain to be a square with length $L=L_{x_1}=L_{x_2}=100\,\text{mm}$. The uncertain diffusion and reaction coefficients are defined as in \eqref{eq:a-def}--\eqref{eq:k-def}, taken to be independent of $t$, i.e., stationary in time. The mean fields $a_0:=\psi_0(\bsx)=0.05$ mm$^2$/day and $\kappa_0:=\xi_0(\bsx)=0.3$ day$^{-1}$ are taken to be in realistic ranges based on previous model calibration studies for gliomas \cite{hormuth2021image}. The fluctuations are prescribed by interleaving two sequences as follow:
\begin{align*}
\psi_{j}(\bsx) &= 
 \begin{cases}
 a_0\,\frac{1}{2}\,k^{-\nu}\sin(\frac{k\,\pi\,x_1}{L_{x_1}})\,\sin(\frac{k\,\pi\,x_2}{L_{x_2}})
 & \mbox{if } j = 2k-1. \\
 0 & \mbox{if } j = 2k,     
 \end{cases} 
 \\[2mm]
 \xi_{j}(\bsx) &= 
 \begin{cases}
 0 & \mbox{if } j = 2k-1, \\   
 \kappa_0\,\frac{1}{2}\,k^{-\nu}\sin(\frac{k\,\pi\,x_1}{L_{x_1}})\, 
 \sin(\frac{k\,\pi\,x_2}{L_{x_2}})
 & \mbox{if } j = 2k,
 \end{cases}
\end{align*}
for $j\geq 1$ with smoothness parameter $\nu=2$.

The initial condition $u_0(\bsx)$ is a Gaussian function in the center of the domain as in \cite{rockne2009mathematical} and the terminal time is taken to be $T=7$ days. Radiotherapy \eqref{eq:f-rt} is applied the first five days and chemotherapy \eqref{eq:f-ct} is administered daily. The radiotherapy dosage is uniformly $z_{\text{rt},k}(\bsx)=2$~Gy and the chemotherapy concentration is assumed to be uniformly $z_{\text{ct},\ell}(\bsx)=1$. The radiosensitivity parameter is $\alpha_{\text{rt}}=0.025\,\text{Gy}^{-1}$ with ratio $\alpha_{\text{rt}}/\beta_{\text{rt}}=10$~Gy following reported values in the literature \cite{rockne2009mathematical}. The chemotherapy efficacy is taken to be $\alpha_{\text{ct}}=0.9$ with clearance rate $\beta_{\text{ct}}=24\, \ln(2)/1.8~\text{day}^{-1}$ for a chemotherapy based on Temozolomide, which has a half-life of 1.8~hours \cite{friedman2000temozolomide}. The implicit Euler method is used to discretize \eqref{eq:pde1} in time with time step $\Delta t = 1/8$ day, such that $\gamma_{\text{rt},\varepsilon}=1/\Delta t$, where $\varepsilon = \Delta t$ matches the resolution of the time-discretization. With these specifications, the time-integrated effect of radiotherapy dose $k$ from \eqref{eq:f-rt} is
\[
 \gamma_{\text{rt},\varepsilon} 
 \int_{\tau_{\text{rt},k}}^{\tau_{\text{rt},k}+\varepsilon}
 [1 - \exp(-\alpha_{\text{rt}} z_{\text{rt},k} - \beta_{\text{rt}} z^2_{\text{rt},k})]\,\rd t
 = 1 - \exp(-\alpha_{\text{rt}} z_{\text{rt},k} - \beta_{\text{rt}} z^2_{\text{rt},k})\approx 0.06,
\]
hence each dose of radiotherapy reduces $u$ by about $6\%$. Similarly, the time-integrated effect of chemotherapy dose $\ell$ from \eqref{eq:f-ct} is
\[
 \alpha_{\text{ct}}
 \int_{\tau_{\text{ct},\ell}}^\infty 
 z_{\text{ct},\ell}
 \exp(-\beta_{\text{ct}}(t-\tau_{\text{ct},\ell})) \,\rd t
 = \frac{\alpha_{\text{ct}}\, z_{\text{ct},\ell}}{\beta_{\text{ct}}}\approx 0.10,
\]
and so each dose of chemotherapy reduces $u$ by about $10\%$, with a slight time delay. Linear finite elements are used to discretize the spatial dimension with mesh size $h=1\,\text{mm}$.

We take the quantity of interest to be the total tumor cellularity at the end of the simulation horizon,
\begin{equation*}
    G(u_s^{\bsy}) := \int_\Omega u_s^{\bsy}(\bsx, T)\,\rd\bsx,
\end{equation*}
with stochastic dimension $s=256$, that is, there are $128$ terms in the expansion for each field.
To obtain an unbiased estimator of $\bbE[G(u_s)] = \int_{U_s} G(u_s^\bsy)\,\rd\bsy$, we average over $R$ i.i.d.~uniform random shifts,
\[
 \overline{Q} := \frac{1}{R} \sum_{r=1}^R Q_N^{(r)} (G(u_s)),
 \qquad 
 Q^{(r)}_{N}[G(u_s)] := \frac{1}{N}\sum^N_{i=1} G\big(u^{\bsy^{(i)}_{\bsDelta_{r}}}\big),
\]
where, similarly to \eqref{eq:qmc-estimate}, $Q^{(r)}_{N}[G(u_s)]$
denotes the unbiased \gls*{qmc} approximation of the integral $\bbE[G(u_s)]$ with the $r$th random shift. Therefore the average $\overline{Q}$ is also an unbiased estimator of $\bbE[G(u_s)]$, and its variance is $1/R$ times the variance for any individual $Q^{(r)}_{N}[G(u_s)]$.
To estimate the root-mean-square error of $\overline{Q}$ we compute the approximation
\begin{equation} \label{eqn:rms}
    \sqrt{\frac{1}{R(R-1)} \sum_{r=1}^R \big\vert \;\overline{Q} - Q_N^{(r)}(G(u_s)) \;\big\vert^2}.
\end{equation}
This quantity should behave like the bound on the root-mean-square error from~\Cref{thm:qmc-error} with an additional $1/\sqrt{R}$ since this is the root-mean-square error of $\overline{Q}$ which uses $R$ shifts instead of a single shift. Here we use $R=16$ random shifts $\bsDelta$ and $N=2^m$ lattice points with $m=[4, 5, \dots, 17]$, for a total of $2^{21}=2,097,152$ \gls*{pde} evaluations. An off-the-shelf embedded lattice rule with weights tailored to a spectral decay rate of $1/j^2$ is utilized to generate the \gls*{qmc} sequence~\cite{cools2006constructing}. For comparison, we use $N^{\text{MC}} = 2^{17}\cdot R = 2^{21}$ quadrature points for \gls*{mc} to match the total number of \gls*{pde} evaluations. Computations are carried out on the Texas Advanced Computing Center's \texttt{Lonestar6} cluster across two compute nodes with 256 total parallel MPI processes, each evaluating one \gls*{pde}. The \gls*{mc} reference calculation with $2^{21}$ evaluations took 1,141 minutes ($\approx 19$~hours), while the QMC shift replicates with $2^{17}$ evaluations averaged 146 minutes ($\approx 2.4$~hours). Code to reproduce the experiments is available at: \hyperlink{https://github.com/gtpash/qmc-tumor}{https://github.com/gtpash/qmc-tumor}. The results are displayed in \Cref{fig:uniform-convergence}. We observe that the QMC root-mean-square error converges at a linear rate $\calO(N^{-1})$ while standard \gls*{mc} converges with the typical $\mathcal{O}(N^{-1/2})$ rate, consistent with the theory.

\begin{figure}[t]
    \centering
    \includegraphics[width=0.5\linewidth]{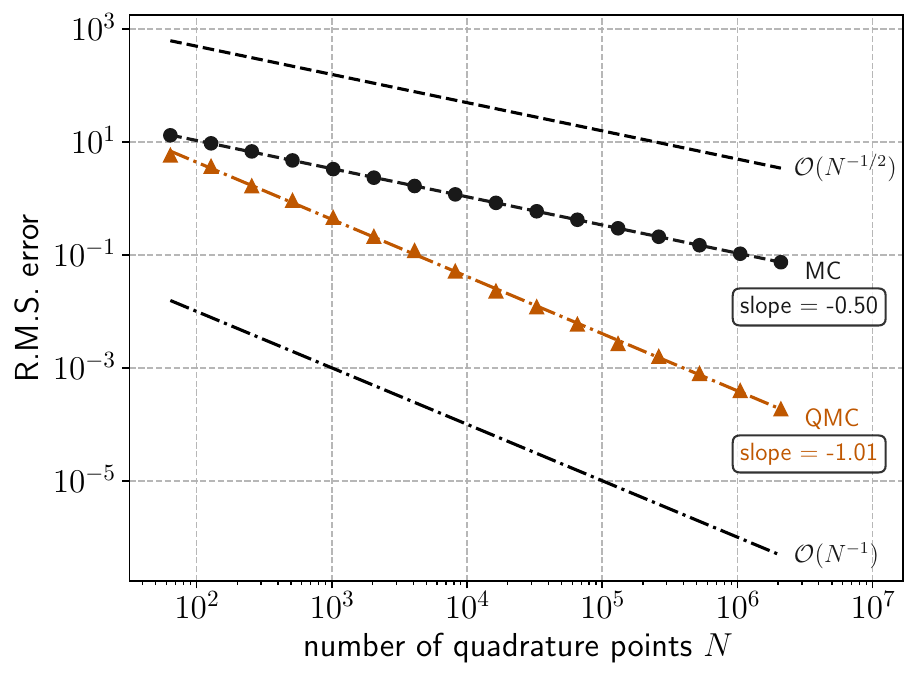}
    \caption{The approximate root-mean-square error \eqref{eqn:rms} of the randomized \gls*{qmc} approximation and the standard error of the \gls*{mc} estimator of the integral $\bbE[G(u_s)]$ for different numbers of total \gls*{pde} evaluations.}
    \label{fig:uniform-convergence}
\end{figure}

\subsection{Experiments with lognormal random fields}
\label{sec:lognormal-numerics}

In this section, we consider a more realistic lognormal model for the uncertain coefficients, despite not yet having theoretical results as in the uniform case. There is QMC theory to handle the lognormal case, see e.g.,~\cite{graham2015quasi,graham2018circulant}, but it requires a new parametric regularity bound for the lognormal model of this problem which is anticipated to be harder and so is left for future work. The following lognormal model has been successfully employed in Bayesian inversion for gliomas in both murine models \cite{liang2023bayesian} and human patients \cite{pash2025predictive},
\begin{align*}
 a^\bsy(\bsx) 
 &= \psi_0(\bsx) + \exp(Z_a(\bsx)),
 \; Z_a(\bsx) = \sum_{j=1}^\infty y_j\,\psi_j(\bsx),
 \; \psi_j = 
 \begin{cases} 
 \sqrt{\mu_{a,k}}\,\varphi_{a,k}(\bsx) & \mbox{if } j = 2k-1, \\
 0 & \mbox{if } j = 2k, 
 \end{cases}
 \\[2mm]
  \kappa^\bsy(\bsx) 
 &= \xi_0(\bsx) + \exp(Z_\kappa(\bsx)), 
 \;\, Z_\kappa(\bsx) = \sum_{j=1}^\infty y_j\,\xi_j(\bsx),
 \;\; \xi_j = 
 \begin{cases} 
 0 & \mbox{if } j = 2k-1, \\
 \sqrt{\mu_{\kappa,k}}\,\varphi_{\kappa,k}(\bsx) & \mbox{if } j = 2k, 
 \end{cases}
\end{align*}
with i.i.d.~$y_{j}\sim\mathcal{N}(0,1)$. Here the parameters are stationary in time once more. The infinite sums are the \gls*{kl} expansions \cite{loeve1978probability2} of the underlying Gaussian random fields $Z_a(\bsx) \sim \mathcal{N}(0,\mathcal{C}_a)$ and $Z_\kappa(\bsx) \sim \mathcal{N}(0,\mathcal{C}_\kappa)$ for the diffusion and proliferation coefficients, with ordered eigenpairs $\{(\mu_{a,k}, \varphi_{a,k})\}_k$ and $\{(\mu_{\kappa,k}, \varphi_{\kappa,k})\}_k$, respectively. We again interleave the two independent sequences. The foundational work of \cite{lindgren2011explicit} established a link between Gaussian random fields with Mat\'ern covariance operators and the solution of a stochastic \gls*{pde}. Following \cite{stuart2010inverse,villa2021hippylib}, the covariance operators $\mathcal{C}_a$ and $\mathcal{C}_\kappa$ are taken to be the inverse of an elliptic operator of the form 
\begin{equation}
\label{eq:covariance}
    \mathcal{C}:=\mathcal{A}^{-\nu}=(-\gamma\,\Delta + \delta\, I)^{-\nu},
\end{equation}
with smoothness parameter $\nu=2$ and hyperparameters $\gamma>0$ and $\delta>0$ that control the correlation length and pointwise variance \cite{villa2024note}. We truncate the expansion and consider the $s$-dimensional subspace spanned by the leading eigenvectors. In particular, we solve a symmetrized (general) eigenvalue problem with a randomized double pass algorithm \cite{villa2021hippylib},
\begin{equation}
\label{eq:gevp}
\mathbf{M}\,\mathbf{\Gamma}\, \boldsymbol{\varphi} = \mu\, \mathbf{M}\, \boldsymbol{\varphi},
\end{equation}
where $\mathbf{\Gamma} := \mathbf{R}^{-1}\mathbf{M}$ is the covariance matrix that arises from the discretization of \eqref{eq:covariance}, with mass matrix $\mathbf{M}_{ij}=\int_\Omega \phi_i(\bsx)\phi_j(\bsx)\,\rd\bsx$ and bilinear form of the elliptic operator given by $\mathbf{R}_{ij}=\int_\Omega \phi_i(\bsx)\mathcal{A}^\nu \phi_j(\bsx)\,\rd\bsx$. Here $\{\phi_i\}_{i=1}^{m}$ is the collection of finite element basis functions. The eigenpairs $\{(\mu_{a,k}, \varphi_{a,k})\}_k$ and $\{(\mu_{\kappa,k}, \varphi_{\kappa,k})\}_k$ obtained by solving \eqref{eq:gevp} may be efficiently computed for generic unstructured domains.

In the numerical experiments, we fix the correlation length to $180$ mm (roughly the largest dimension of a human brain). The pointwise variance is set to $0.2336$ and $0.0682$ for the diffusion and proliferation coefficients, respectively. As in \Cref{sec:uniform-numerics}, we take the mean fields to be $a_0:=\psi_0(\bsx)=0.05$ mm$^2$/day and $\kappa_0:=\xi_0(\bsx)=0.3$ day$^{-1}$. The domain is taken to be a two-dimensional slice of a human left-hemisphere \cite{mardal2022mathematical}. The computational mesh is visualized along with a representative sample in Figure~\ref{fig:lognormal-samples}. The initial condition $u_0(\bsx)$ is once more a Gaussian function centered at the geometric center of the domain with amplitude $0.8$ and width $5.0$ mm. The treatment schedule and radiosensitivity parameters are as in \Cref{sec:uniform-numerics}. A time step of $\Delta t=1/8$ day is used for the temporal discretization and the unstructured mesh contains $9,995$ vertices.

We empirically investigate the convergence rate using $R=16$ random shifts $\bsDelta$ and $N=2^m$ lattice points with $m=[6, 7, \dots 17]$. Each \gls*{kl} expansion is truncated at $512$ terms, for a stochastic dimension of $s=1024$. In the lognormal case, the shifted QMC points are 
$\bsy^{(i)}_{\bsDelta} 
    := \boldsymbol{\Phi}^{-1}\left((\tfrac{i}{N}\,\bsz + \bsDelta) \bmod 1\right)$,
where $\boldsymbol{\Phi}^{-1}(\cdot)$ denotes the inverse cumulative normal distribution applied componentwise to a vector. The same embedded lattice rule as in the uniform case is used to construct the \gls*{qmc} sequence. Computations are carried out in the same fashion as in \Cref{sec:uniform-numerics}. The reference calculation took 1,417 minutes (23.6~hours) and the QMC evaluations averaged 165 minutes (2.8~hours). We attribute the slight increase in computational cost to the larger stochastic dimension which resolves finer features as well as the non-uniform domain. The empirical spectrum of the covariance operators $\mathcal{C}_a$ and $\mathcal{C}_\kappa$ and convergence results are shown in Figure~\ref{fig:lognormal-results}. We observe that the QMC root-mean-square error converges at a faster rate of $\mathcal{O}(N^{-0.86})$ than the MC rate of $\mathcal{O}(N^{-1/2})$. The square-root of the eigenvalues decay like $k^{-1.17}$ and so we expect the QMC convergence rate to be $1/p-1/2 \approx 1.17-1/2 = 0.67$, if the theory for uniform fields also holds for lognormal fields (as proven for some elliptic linear problems). Hence our results are very encouraging and prompt further analysis. 

\begin{figure}[t]
    \centering
    \includegraphics[width=\linewidth]{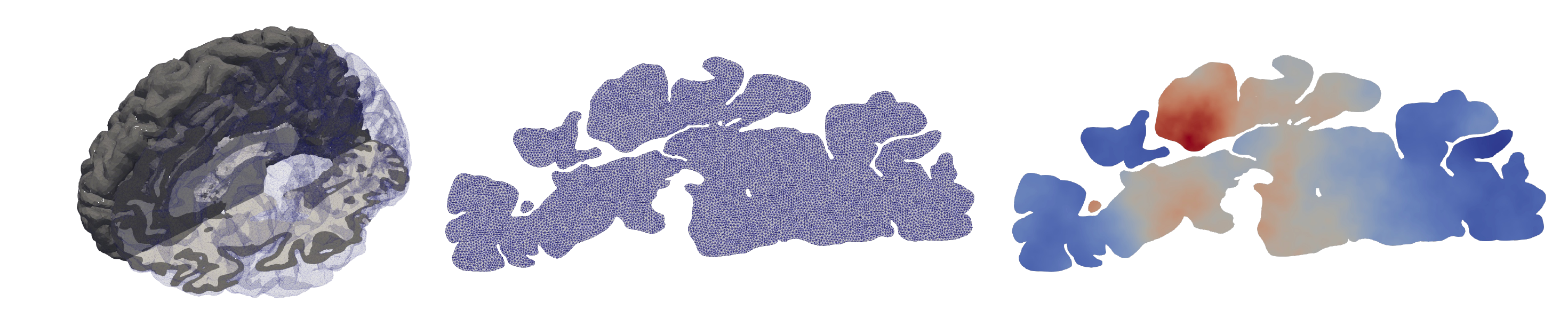}
    \caption{Left: Visualization of anatomy and slice of the left-hemisphere. Center: The computational mesh used for simulations, with $9,995$ vertices. Right: A representative sample of the diffusion coefficient.}
    \label{fig:lognormal-samples}
\end{figure}

\begin{figure}[t]
     \centering
     \begin{subfigure}{0.49\textwidth}
         \centering
         \includegraphics[width=\linewidth]{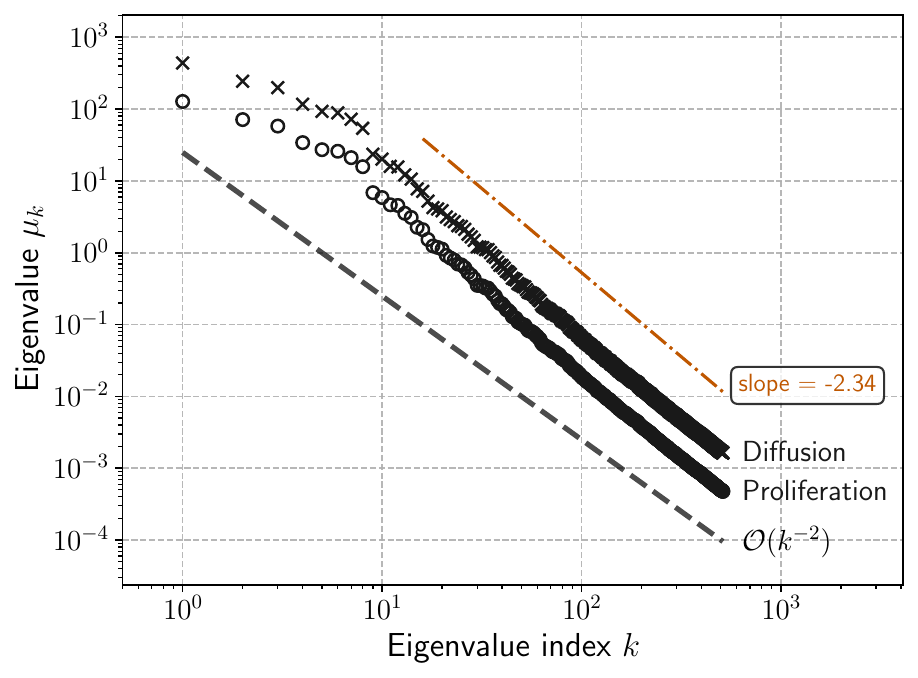}
     \end{subfigure}
     \begin{subfigure}{0.49\textwidth}
         \centering
         \includegraphics[width=\linewidth]{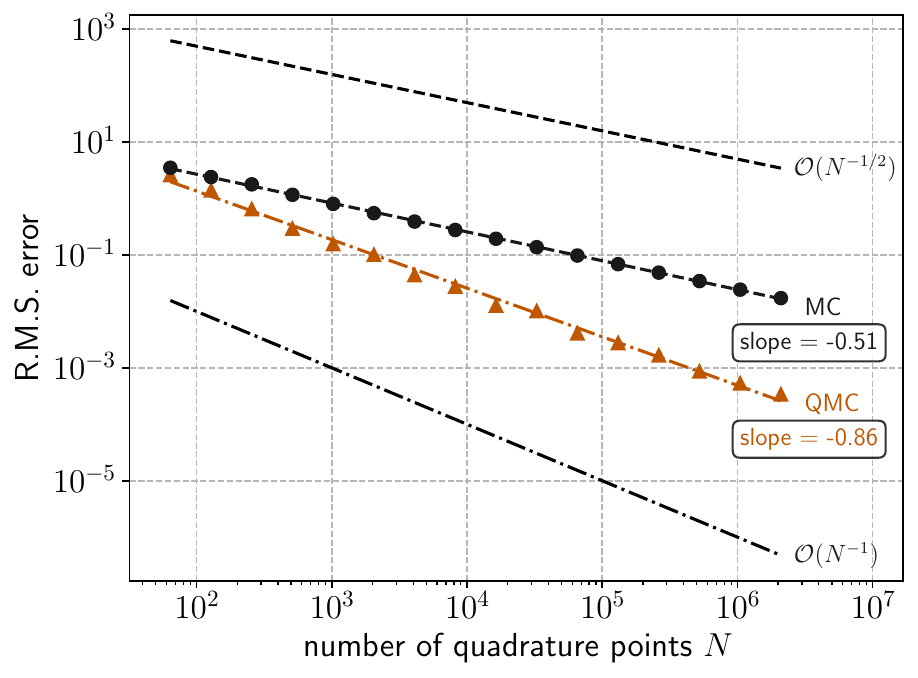}
     \end{subfigure}
     \caption{Left: The empirical spectral decay of the two Gaussian random fields compared to the theoretical rate accompanied. Right: The approximate root-mean-square error \eqref{eqn:rms} for the \gls*{qmc} approximation and the standard error for the \gls*{mc} estimator of the integral $\bbE[G(u_s)]$.}
     \label{fig:lognormal-results}
\end{figure}

\section{Conclusions}
\label{sec:conclusions}
In this work, we developed a \gls*{qmc} method for a semi-linear parabolic \gls*{pde} model of tumor growth. The method was applied to the high-dimensional integral arising from estimating the expected value of a smooth \gls*{qoi}. We developed a well-posedness argument for this particular class of problems and derived an \textit{a priori} bound with an explicit constant. Moreover, in the case of uniform random fields we obtained a parametric regularity bound to provide theoretical error bounds for the \gls*{qmc} approximation. The convergence rate was then numerically verified, in particular, we observe linear convergence. Additionally, we numerically studied application to lognormal random fields and once more observed faster than MC convergence. Future work will focus on extending approximation theory to the lognormal case and application to an optimal control problem. Furthermore, a rigorous analysis of the dimension truncation error as well as the numerical discretization of the \gls*{pde} should be executed to assess their contributions to the overall discretization error.

\section*{Acknowledgments}
The authors give heartfelt thanks to both referees for their detailed suggestions to improve the paper, especially to one referee for explaining how to prove that the unique weak solution is always nonnegative. We are also grateful for insightful discussions with Albert Cohen, Wolfgang Dahmen, Mike Giles, Maria Pia Gualdani, Pierre L'Ecuyer, Dingcheng Luo, Fabio Nobile, Mark Peletier, William Porteous, and Christian Wieners, as well as other participants at the 2025 Oberwolfach Workshop on Computation and Learning in High Dimensions.

\bibliographystyle{siamplain}
\bibliography{references}

\begin{thebibliography}{10}

\bibitem{alfonso2017biology}
{\sc J.~Alfonso, K.~Talkenberger, M.~Seifert, B.~Klink, A.~Hawkins-Daarud, K.~Swanson, H.~Hatzikirou, and A.~Deutsch}, {\em The biology and mathematical modelling of glioma invasion: a review}, Journal of the Royal Society Interface, 14 (2017), p.~20170490.

\bibitem{an2025sparse}
{\sc X.~An, J.~Dick, M.~Feischl, A.~Scaglioni, and T.~Tran}, {\em {Sparse Grid Approximation of Nonlinear SPDEs: The Landau--Lifshitz--Gilbert Equation}}, SIAM/ASA J. Uncertain. Quantif., 13 (2025), pp.~472--517.

\bibitem{borasi2016modelling}
{\sc G.~Borasi and A.~Nahum}, {\em Modelling the radiotherapy effect in the reaction-diffusion equation}, Physica Medica, 32 (2016), pp.~1175--1179.

\bibitem{CCMR25}
{\sc G.~Cavalleri, P.~Colli, A.~Miranville, and E.~Rocca}, {\em On a brain tumor growth model with lactate metabolism, viscoelastic effects, and tissue damage}, Nonlinear Anal. Real World Appl., 87 (2026), p.~104419.

\bibitem{chaplain1996avascular}
{\sc M.~Chaplain}, {\em Avascular growth, angiogenesis and vascular growth in solid tumours: The mathematical modelling of the stages of tumour development}, Math. Comput. Model., 23 (1996), pp.~47--87.

\bibitem{CL24a}
{\sc A.~Chernov and T.~L{\^e}}, {\em Analytic and {G}evrey class regularity for parametric elliptic eigenvalue problems and applications}, SIAM J. Numer. Anal., 62 (2024), pp.~1874--1900.

\bibitem{CL24b}
{\sc A.~Chernov and T.~L{\^e}}, {\em Analytic and {G}evrey class regularity for parametric semilinear reaction-diffusion problems and applications in uncertainty quantification}, Comput. Math. Appl., 164 (2024), pp.~116--130.

\bibitem{chernov2025gevrey}
{\sc A.~Chernov and T.~L\^e}, {\em Gevrey class regularity for steady-state incompressible {N}avier-{S}tokes equations in parametric domains and related models}, arXiv preprint arXiv:2504.13753,  (2025).

\bibitem{clatz2005realistic}
{\sc O.~Clatz, M.~Sermesant, P.-Y. Bondiau, H.~Delingette, S.~K. Warfield, G.~Malandain, and N.~Ayache}, {\em Realistic simulation of the 3-{D} growth of brain tumors in {MR} images coupling diffusion with biomechanical deformation}, IEEE Trans. Med. Imaging, 24 (2005), pp.~1334--1346.

\bibitem{CGLMRR20}
{\sc P.~Colli, H.~Gomez, G.~Lorenzo, G.~Marinoschi, A.~Reali, and E.~Rocca}, {\em Mathematical analysis and simulation study of a phase-field model of prostate cancer growth with chemotherapy and antiangiogenic therapy effects}, Math. Models Methods Appl. Sci., 30 (2020), pp.~1253--1295.

\bibitem{cools2006constructing}
{\sc R.~Cools, F.~Y. Kuo, and D.~Nuyens}, {\em Constructing embedded lattice rules for multivariate integration}, SIAM J. Sci. Comput., 28 (2006), pp.~2162--2188.

\bibitem{DL92}
{\sc R.~Dautry and J.-L. Lions}, {\em {Mathematical analysis and numerical methods for science and technology: Evolution problems I}}, Springer-Verlag, 1992.

\bibitem{dick2014higher}
{\sc J.~Dick, F.~Y. Kuo, Q.~T. Le~Gia, D.~Nuyens, and C.~Schwab}, {\em Higher order {QMC} {P}etrov--{G}alerkin discretization for affine parametric operator equations with random field inputs}, SIAM J. Numer. Anal., 52 (2014), pp.~2676--2702.

\bibitem{Eva10}
{\sc L.~C. Evans}, {\em Partial Differential Equations}, vol.~19, American Mathematical Society, 2022.

\bibitem{fisher1937wave}
{\sc R.~A. Fisher}, {\em The wave of advance of advantageous genes}, Annals of Eugenics, 7 (1937), pp.~355--369.

\bibitem{friedman2000temozolomide}
{\sc H.~S. Friedman, T.~Kerby, and H.~Calvert}, {\em Temozolomide and treatment of malignant glioma}, Clinical Cancer Research, 6 (2000), pp.~2585--2597.

\bibitem{gantner2018}
{\sc R.~N. Gantner, L.~Herrmann, and C.~Schwab}, {\em {Quasi--Monte Carlo Integration for Affine-Parametric, Elliptic PDEs: Local Supports and Product Weights}}, SIAM J. Numer. Anal., 56 (2018), pp.~111--135.

\bibitem{GT77}
{\sc D.~Gilbarg and N.~S. Trudinger}, {\em {Elliptic Partial Differential Equations of Second Order}}, vol.~224, Springer, 1977.

\bibitem{gilbert2019evp}
{\sc A.~D. Gilbert, I.~G. Graham, F.~Y. Kuo, I.~H. Sloan, and R.~Scheichl}, {\em Analysis of quasi-{M}onte {C}arlo methods for elliptic eigenvalue problems with stochastic coefficients}, Numer. Math., 142 (2019), pp.~863--915.

\bibitem{graham2015quasi}
{\sc I.~G. Graham, F.~Y. Kuo, J.~A. Nichols, R.~Scheichl, C.~Schwab, and I.~H. Sloan}, {\em Quasi-{M}onte {C}arlo finite element methods for elliptic {PDE}s with lognormal random coefficients}, Numer. Math., 131 (2015), pp.~329--368.

\bibitem{graham2018circulant}
{\sc I.~G. Graham, F.~Y. Kuo, D.~Nuyens, R.~Scheichl, and I.~H. Sloan}, {\em Circulant embedding with {QMC}: analysis for elliptic {PDE} with lognormal coefficients}, Numer. Math., 140 (2018), pp.~479--511.

\bibitem{guth2021quasi}
{\sc P.~A. Guth, V.~Kaarnioja, F.~Y. Kuo, C.~Schillings, and I.~H. Sloan}, {\em A quasi-{M}onte {C}arlo method for optimal control under uncertainty}, SIAM/ASA J. Uncertain. Quantif., 9 (2021), pp.~354--383.

\bibitem{GKKSS24}
{\sc P.~A. Guth, V.~Kaarnioja, F.~Y. Kuo, C.~Schillings, and I.~H. Sloan}, {\em Parabolic {PDE}-constrained optimal control under uncertainty with entropic risk measure using quasi-{M}onte {C}arlo integration}, Numer. Math., 156 (2024), pp.~565--608.

\bibitem{hakula2024}
{\sc H.~Hakula, H.~Harbrecht, V.~Kaarnioja, F.~Y. Kuo, and I.~H. Sloan}, {\em {Uncertainty quantification for random domains using periodic random variables}}, Numer. Math., 156 (2024), pp.~273--317.

\bibitem{harbrecht2024}
{\sc H.~Harbrecht, M.~Schmidlin, and C.~Schwab}, {\em {The Gevrey class implicit mapping theorem with application to UQ of semilinear elliptic PDEs}}, Math. Models Methods Appl. Sci., 34 (2024), pp.~881--917.

\bibitem{he2024analytic}
{\sc Y.~He and C.~Schwab}, {\em Analytic regularity and solution approximation for a semilinear elliptic partial differential equation in a polygon}, Calcolo, 61 (2024), p.~11.

\bibitem{hormuth2021image}
{\sc D.~A. Hormuth, K.~A. Al~Feghali, A.~M. Elliott, T.~E. Yankeelov, and C.~Chung}, {\em Image-based personalization of computational models for predicting response of high-grade glioma to chemoradiation}, Scientific Reports, 11 (2021), p.~8520.

\bibitem{hormuth2015predicting}
{\sc D.~A. Hormuth~II, J.~A. Weis, S.~L. Barnes, M.~I. Miga, E.~C. Rericha, V.~Quaranta, and T.~E. Yankeelov}, {\em Predicting in vivo glioma growth with the reaction diffusion equation constrained by quantitative magnetic resonance imaging data}, Physical Biology, 12 (2015), p.~046006.

\bibitem{jarrett2021quantitative}
{\sc A.~M. Jarrett, A.~S. Kazerouni, C.~Wu, J.~Virostko, A.~G. Sorace, J.~C. DiCarlo, D.~A. Hormuth, D.~A. Ekrut, D.~Patt, B.~Goodgame, S.~Avery, and T.~E. Yankeelov}, {\em Quantitative magnetic resonance imaging and tumor forecasting of breast cancer patients in the community setting}, Nature Protocols, 16 (2021), pp.~5309--5338.

\bibitem{kolmogoroff1988study}
{\sc A.~Kolmogoroff, I.~Petrovsky, and N.~Piscounoff}, {\em Study of the diffusion equation with growth of the quantity of matter and its application to a biology problem}, in Dynamics of Curved Fronts, Elsevier, 1988, pp.~105--130.

\bibitem{kuo2016application}
{\sc F.~Y. Kuo and D.~Nuyens}, {\em Application of quasi-{M}onte {C}arlo methods to elliptic pdes with random diffusion coefficients: a survey of analysis and implementation}, Found. Comput. Math., 16 (2016), pp.~1631--1696.

\bibitem{kuo2012quasi}
{\sc F.~Y. Kuo, C.~Schwab, and I.~H. Sloan}, {\em Quasi-{M}onte {C}arlo finite element methods for a class of elliptic partial differential equations with random coefficients}, SIAM J. Numer. Anal., 50 (2012), pp.~3351--3374.

\bibitem{liang2023bayesian}
{\sc B.~Liang, J.~Tan, L.~Lozenski, D.~A. Hormuth, T.~E. Yankeelov, U.~Villa, and D.~Faghihi}, {\em Bayesian inference of tissue heterogeneity for individualized prediction of glioma growth}, IEEE Trans. Med. Imaging, 42 (2023), pp.~2865--2875.

\bibitem{lindgren2011explicit}
{\sc F.~Lindgren, H.~Rue, and J.~Lindstr{\"o}m}, {\em An explicit link between {G}aussian fields and {G}aussian {M}arkov random fields: the stochastic partial differential equation approach}, J. R. Stat. Soc. Ser. B. Stat. Methodol., 73 (2011), pp.~423--498.

\bibitem{loeve1978probability2}
{\sc M.~Lo{\`e}ve}, {\em Probability Theory II}, vol.~46 of Graduate Texts in Mathematics, Springer, New York, 4~ed., 1978.

\bibitem{lorenzo2019computer}
{\sc G.~Lorenzo, T.~J. Hughes, P.~Dominguez-Frojan, A.~Reali, and H.~Gomez}, {\em Computer simulations suggest that prostate enlargement due to benign prostatic hyperplasia mechanically impedes prostate cancer growth}, Proc. Nat. Acad. Sci., 116 (2019), pp.~1152--1161.

\bibitem{mardal2022mathematical}
{\sc K.-A. Mardal, M.~E. Rognes, T.~B. Thompson, and L.~M. Valnes}, {\em {Mathematical Modeling of the Human Brain: From Magnetic Resonance Images to Finite Element Simulation}}, vol.~10, Springer Nature, 2022.

\bibitem{mcmahon2018linear}
{\sc S.~J. McMahon}, {\em The linear quadratic model: usage, interpretation and challenges}, Physics in Medicine \& Biology, 64 (2018), p.~01TR01.

\bibitem{murray2003mathematical}
{\sc J.~D. Murray}, {\em Mathematical {B}iology {II}: {S}patial {M}odels and {B}iomedical {A}pplications}, vol.~18, Springer Science \& Business Media, 2003.

\bibitem{nguyen2024evp}
{\sc V.~K. Nguyen}, {\em Analyticity of parametric elliptic eigenvalue problems and applications to quasi-{M}onte {C}arlo methods}, Complex Var. Elliptic Equ., 69 (2024), pp.~1--21.

\bibitem{pash2025predictive}
{\sc G.~Pash, U.~Villa, D.~A. Hormuth~II, T.~E. Yankeelov, and K.~Willcox}, {\em Predictive digital twins with quantified uncertainty for patient-specific decision making in oncology}, arXiv preprint arXiv:2505.08927,  (2025).

\bibitem{Per15}
{\sc B.~Perthame}, {\em Parabolic Equations in Biology}, Springer, 2015.

\bibitem{quarteroni1994numerical}
{\sc A.~Quarteroni and A.~Valli}, {\em Numerical Approximation of Partial Differential Equations}, Springer, 1994.

\bibitem{rockne2009mathematical}
{\sc R.~Rockne, E.~Alvord~Jr, J.~K. Rockhill, and K.~R. Swanson}, {\em A mathematical model for brain tumor response to radiation therapy}, J. Math. Biol., 58 (2009), pp.~561--578.

\bibitem{SchSt09}
{\sc C.~Schwab and R.~Stevenson}, {\em {Space-time adaptive wavelet methods for parabolic evolution problems}}, Math.~Comp., 78 (2009), pp.~1293--1318.

\bibitem{stuart2010inverse}
{\sc A.~M. Stuart}, {\em Inverse problems: {A} {B}ayesian perspective}, Acta Numer., 19 (2010), pp.~451--559.

\bibitem{stupp2005radiotherapy}
{\sc R.~Stupp, W.~P. Mason, M.~J. Van Den~Bent, M.~Weller, B.~Fisher, M.~J. Taphoorn, K.~Belanger, A.~A. Brandes, C.~Marosi, U.~Bogdahn, et~al.}, {\em Radiotherapy plus concomitant and adjuvant temozolomide for glioblastoma}, New England Journal of Medicine, 352 (2005), pp.~987--996.

\bibitem{swanson2000quantitative}
{\sc K.~R. Swanson, E.~C. Alvord~Jr, and J.~D. Murray}, {\em A quantitative model for differential motility of gliomas in grey and white matter}, Cell Proliferation, 33 (2000), pp.~317--329.

\bibitem{Tro10}
{\sc F.~Tr{\"o}ltzsch}, {\em Optimal Control of Partial Differential Equations: Theory, Methods, and Applications}, vol.~112, American Mathematical Society, 2010.

\bibitem{villa2024note}
{\sc U.~Villa and T.~O'Leary-Roseberry}, {\em A note on the relationship between {PDE}-based precision operators and {M}at\'ern covariances}, arXiv preprint arXiv:2407.00471,  (2024).

\bibitem{villa2021hippylib}
{\sc U.~Villa, N.~Petra, and O.~Ghattas}, {\em {HIPPYlib}: an extensible software framework for large-scale inverse problems governed by {PDEs}: {P}art {I}: deterministic inversion and linearized {B}ayesian inference}, ACM Trans. Math. Software, 47 (2021), pp.~1--34.

\end{thebibliography}

\end{document}